\documentclass[12pt,
reqno]{amsart}

\usepackage[cp1251]{inputenc}
\usepackage{amsmath,amsxtra,amssymb,eufrak}
\usepackage[all]{xy}
\usepackage[active]{srcltx} 
\usepackage{mathrsfs}
\usepackage[russian]{babel}

\theoremstyle{plain}
\newtheorem{thm}{Теорема}[section]
\newtheorem{lm}[thm]{Лемма}
\newtheorem{co}[thm]{Следствие}
\newtheorem{pr}[thm]{Предложение}

\theoremstyle{definition}

\newtheorem{df}[thm]{Определение}
\newtheorem{exm}[thm]{Пример}
\newtheorem{rem}[thm]{Замечание}
\numberwithin{equation}{section}

\title{Разложение алгебры аналитических функционалов на связной комплексной группе Ли и е\"{е} пополнений в итерированные аналитические смэш-произведения}

\author{О. Ю. Аристов}
\address{Institute for Advanced Study in Mathematics, Harbin Institute of Technology, Harbin 150001, China}
\email{aristovoyu@inbox.ru}
\keywords{Аналитическое смэш-произведение, топологическая алгебра Хопфа, комплексная группа Ли, субмультипликативный вес с экспоненциальным искривлением, функция длины, аналитический функционал, оболочка Аренса-Майкла, оболочка относительно  класса банаховых PI-алгебр}
\begin{document}

 \begin{abstract}
Показано, что  разложение комплексной группы Ли $G$ в полупрямое произведение порождает разложение её алгебры аналитических функционалов ${\mathscr A}(G)$ в аналитическое смэш-произведение в смысле Пирковского. Кроме того, даны достаточные условия того, что полупрямое произведение порождает аналогичные разложения некоторых алгебр Аренса-Майкла, которые являются пополнениями ${\mathscr A}(G)$. Основной результат: если $G$ связна, то её линеаризация допускает разложение в итерированное полупрямое произведение (соответствующий композиционный ряд содержит абелевы и полупростой факторы), которое индуцирует разложение в итерированное аналитическое смэш-произведение алгебр из некоторого класса пополнений ${\mathscr A}(G)$. Рассматривая крайние случаи,  оболочку ${\mathscr A}(G)$ относительно  класса всех  банаховых алгебр (она же оболочка Аренса-Майкла) и оболочку относительно  класса банаховых PI-алгебр (новое понятие, которое вводится в этой статье), мы получаем, в частности, их разложения в итерированные аналитические смэш-произведения.
 \end{abstract}

 \maketitle
 \markright{Разложение алгебры аналитических функционалов}

\section*{Введение}

Хорошо известно, что разложение группы в полупрямое произведение индуцирует разложение её групповой алгебры в смэш-произведение (см., например, \cite[с.\,351]{Ra}). В свою очередь, разложение алгебры Ли в полупрямое произведение индуцирует разложение её универсальной обертывающей алгебры в смэш-произведение (см., например,  \cite[с.\,33--34, 1.7.11(iv)]{MR87}). В этой статье получены аналогичные результаты для групповых алгебр в контексте комплексного анализа, с акцентом на возможность разложения в итерированное аналитическое смэш-произведение с достаточно просто устроенными факторами.

Мы рассмотрим три локально выпуклые алгебры, ассоциированные с комплексной группой Ли $G$: алгебру аналитических функционалов ${\mathscr A}(G)$, а также два её пополнения --- алгебры Аренса-Майкла $\widehat{{\mathscr A}}(G)$ и $\widehat{{\mathscr A}}(G)^{\mathrm{PI}}$ (которые являются частными случаями пополнений более общего вида, см. ниже).

(1)~В случае $G=\mathbb{C}^n$ алгебра ${\mathscr A}(G)$ известна очень давно, а для произвольной комплексной группы Ли, введена Литвиновым более 50-ти лет назад, см.~\cite{Li70,Li72}, и с тех получила некоторое (хотя и недостаточное по сравнению с алгебрами распределений) внимание математиков.

(2)~Её оболочка Аренса-Майкла $\widehat{{\mathscr A}}(G)$ (т.е. пополнение относительно семейства всех  субмультипликативных непрерывных преднорм) интересна, в частности, тем, что связана с конструкцией голоморфной рефлексивности, см. \cite{Ak08}, а также~\cite{ArHR}.  Для связной $G$ она рассмотрена автором в \cite{ArAnF} и \cite{AHHFG} с точки зрения large-scale геометрии и функционального анализа.

(3)~Алгебра $\widehat{{\mathscr A}}(G)^{\mathrm{PI}}$, которая по определению есть пополнение ${\mathscr A}(G)$
относительно семейства тех субмультипликативных непрерывных преднорм, пополнения по которым являются PI-алгебрами (удовлетворяют полиномиальному тождеству), по-видимому, ранее не исследовалась.
Как $\widehat{{\mathscr A}}(G)^{\mathrm{PI}}$, так и аналогичные оболочки других алгебр относительно класса банаховых PI-алгебр (см. определение~\ref{BPIun}) интересны как объекты некоммутативной комплексно-аналитической геометрии. Представляется, что таких оболочек изучение может пролить свет на взаимосвязи разных её разделов.

С точки зрения теории представлений указанные алгебры характеризуются следующими свойствами. Голоморфным представлениям~$G$ в локально выпуклых пространствах соответствуют непрерывные представления ${\mathscr A}(G)$ \cite[предложение~5]{Li70}, а в банаховых --- непрерывные представления $\widehat{{\mathscr A}}(G)$. Если в последнем случае дополнительно потребовать, чтобы образ каждого представления  удовлетворял полиномиальному тождеству, то мы получим $\widehat{{\mathscr A}}(G)^{\mathrm{PI}}$.

Далее мы предполагаем, что $G$ связна и линейна. Как показано в \S\,\ref{s:UNAL}, $\widehat{{\mathscr A}}(G)$ и $\widehat{{\mathscr A}}(G)^{\mathrm{PI}}$ являются крайними случаями общего семейства алгебр Аренса-Майкла, обладающих некоторыми универсальными свойствами. Члены этого семейства обозначаются через ${\mathscr A}_{\omega_{max}^\infty}(G)$,  где $\omega_{max}$ ---  субмультипликативный вес, зависящий от выбора нильпотентной подгруппы, промежуточной между экспоненциальным и нильпотентным радикалами, а именно, $\omega_{max}$ есть максимальный субмультипликативный вес с экспоненциальным искривлением на этой подгруппе. Так же как  $\widehat{{\mathscr A}}(G)$ и $\widehat{{\mathscr A}}(G)^{\mathrm{PI}}$, алгебры вида ${\mathscr A}_{\omega_{max}^\infty}(G)$ могут быть рассмотрены с точки зрения large-scale геометрии. (Отметим, что тесная связь между банаховыми алгебрами, являющимися пополнениями ${\mathscr A}(G)$,  и субмультипликативными весами была обнаружена Акбаровым, который использовал термин ``полухарактер'',  см.~\cite{Ak08})

Здесь предложен единый подход к описанию структуры всех этих алгебр. Он основан на аналитической версии конструкции смэш-произведения, предложенной Пирковским в~\cite{Pi4}. Наша основная цель --- получить разложения в итерированное аналитическое смэш-произведение с  факторами  простейшего вида. Мотивировкой здесь являются вопросы, связанные с гомологическими эпиморфизмами. (Автор планирует посвятить этой теме отдельную статью.) Структурная теория связных линейных групп обеспечивает  разложение~$G$  в итерированное полупрямое произведение абелевых и полупростого факторов, см., например,  \cite[\S\,16.3]{HiNe}. Поэтому задача сводится к выяснению условий, при которых полупрямое произведение порождает аналитическое смэш-произведение.

В этой статье существенно используются результаты автора, полученные в \cite{ArPiLie,ArUAPI}. В \cite{ArPiLie} установлена связь с теорией PI-алгебр, а в \cite{ArUAPI} получено асимптотическое разложение максимальной функции длины с экспоненциальным искривлением --- утверждение, которое является основным техническим средством в наших рассуждениях.

Кратко о содержании параграфов. В~\S\,\ref{s:ASP} включены необходимые сведения об аналитических смэш-произведениях.
В~\S\,\ref{s:DAAF} доказано, что разложение комплексной группы Ли в полупрямое произведение порождает
разложение алгебры аналитических функционалов в аналитическое смэш-произведение.
В~\S\,\ref{s:DCASW} даны достаточные условия, при которых полупрямое произведение порождает разложение пополнения алгебры аналитических функционалов, ассоциированного с субмультипликативным весом.
Базовый технический результат, теорема о разложении веса c экспоненциальным искривлением на промежуточной нильпотентной подгруппе, содержится в~\S\,\ref{s:TDWED}. В~\S\,\ref{s:UNAL} обсуждаются универсальные алгебры, в частности, алгебры вида $\widehat{{\mathscr A}}(G)^{\mathrm{PI}}$. В~\S\,\ref{s:DISPAM} доказаны основные утверждения статьи: о разложении в итерированное смэш-произведение алгебр Аренса-Майкла для алгебры ${\mathscr A}_{\omega_{max}^\infty}(G)$, где $\omega_{max}$ ---  максимальный вес с экспоненциальным искривлением на промежуточной нильпотентной подгруппе, а также указан вид, который оно принимает в частных случаях для $\widehat{{\mathscr A}}(G)$ и $\widehat{{\mathscr A}}(G)^{\mathrm{PI}}$.

Автор благодарен А.\,Ю.~Пирковскому, разъяснившему некоторые вопросы, связанные с определением смэш-произведения в~\cite{Pi3}. Также автор хотел бы выразить признательность рецензенту за высококачественный отзыв, в котором были указаны многочисленные неточности и некорректные высказывания. В частности, замечания рецензента мотивировали автора использовать в статье обобщенные обозначения Свидлера и ввести понятие ``асимптотически симметричный вес''.

\section{Аналитические смэш-произведения}
\label{s:ASP}

В этом разделе кратко опиcана аналитическая версия конструкции смэш-про\-из\-ведения для интересующего нас класса топологических алгебр Хопфа. Она была впервые рассмотрена Пирковским в~\cite{Pi4}.

Для нас основными объектами являются алгебры, коалгебры, биалгебры и алгебры Хопфа (и модули над ними) в симметрической моноидальной категории полных локально выпуклых пространств над $\mathbb{C}$ с бифунктором $(-)\mathbin{\widehat{\otimes}} (-)$  полного проективного тензорного произведения. Будем называть их соответственно $\mathbin{\widehat{\otimes}}$-алгебрами, $\mathbin{\widehat{\otimes}}$-коалгебрами, $\mathbin{\widehat{\otimes}}$-биалгебрами и $\mathbin{\widehat{\otimes}}$-алгебрами Хопфа (при чтении вслух символ $\mathbin{\widehat{\otimes}}$ обычно заменяется на слово ``топологическая''). В частности, $\mathbin{\widehat{\otimes}}$-алгебра --- это полная локально выпуклая алгебра с совместно непрерывным умножением. Все ассоциативные алгебры предполагаются имеющими единицу. Мы обозначаем умножение (точнее, его линеаризацию), коумножение, коединицу и антипод через $\mu$, $\Delta$, $\varepsilon$ и $S$, при необходимости с добавлением индекса, указывающего на конкретный объект. Все гомоморфизмы в категориях $\mathbin{\widehat{\otimes}}$-алгебр и т.д. (также как и $\mathbin{\widehat{\otimes}}$-модулей над ними) предполагаются непрерывными и сохраняющими единицу.

\subsection*{Обобщенные обозначения Свидлера}
Тензорная запись операции коумножения иногда бывает очень громоздкой. Поэтому в теории коалгебр и, в частности, алгебр Хопфа, широко используются обозначения Свидлера.  Как отмечено в \cite[\S\,2.4]{Ak08}, обозначения Свидлера могут быть обобщены на случай топологических алгебр Хопфа. Ниже описан вариант, пригодный для наших целей.

В классической нотации Свидлера коумножение на коалгебре имеет вид
\begin{equation}\label{SwDe}
\Delta(h)=\sum_{(h)}h_{(1)}\otimes h_{(2)}.
\end{equation}
Индекс $(h)$  под суммой может быть опущен, как и сам знак суммы.

Мы используем этот тип обозначений также для $\mathbin{\widehat{\otimes}}$-коалгебры.
Если она является пространством Фреше, то мы можем понимать~\eqref{SwDe} как сходящийся ряд. Однако в общем случае элемент проективного тензорного произведения не обязательно представляется в виде ряда.
Поэтому запишем его как предел направленности конечных сумм элементарных тензоров:
$$
\Delta(h)=\lim_\nu \sum_{i=1}^{n_\nu} h_{(1)}^{\nu,i}\otimes h_{(2)}^{\nu,i}.
$$
В случае произвольной $\mathbin{\widehat{\otimes}}$-коалгебры  формула в~\eqref{SwDe} является обозначением этого предела.

Например, аксиому антипода для $\mathbin{\widehat{\otimes}}$-алгебры Хопфа можно записать в том же виде, что и для обычной алгебры Хопфа:
$$
\sum S(h_{(1)})h_{(2)}=\varepsilon(h)1=\sum h_{(1)}S(h_{(2)}).
$$
Она расшифровывается как
$$
\lim_\nu \sum_{i=1}^{n_\nu} S(h_{(1)}^{\nu,i})h_{(2)}^{\nu,i}=\varepsilon(h)1=\lim_\nu \sum_{i=1}^{n_\nu} h_{(1)}^{\nu,i}S(h_{(2)}^{\nu,i}).
$$

При использовании итераций коумножения обобщённые обозначения Свидлера позволяют упростить ещё более громоздкие формулы.
Например, аксиома коаcсоциативности, т.е. равенство $(1\otimes \Delta)\Delta(h)=(\Delta\otimes 1)\Delta(h)$, в этих обозначениях имеет вид
$$
\sum h_{(1)}\otimes \left(\sum(h_{(2)})_{(1)}\otimes (h_{(2)})_{(2)}\right)= \left(\sum(h_{(1)})_{(1)}\otimes (h_{(1)})_{(2)}\right)\otimes h_{(2)}.
$$
Детальная расшифровка выглядит так:
\begin{multline}\label{detcoass}
\lim_\nu \left(\sum_{i=1}^{n_\nu} h_{(1)}^{\nu,i}\otimes \left(\lim_\mu \sum_{j=1}^{m_\mu} (h_{(2)}^{\nu,i})_{(1)}^{\mu,j}\otimes (h_{(2)}^{\nu,i})_{(2)}^{\mu,j}\right)\right)=\\
=
\lim_\nu \left(\left(\lim_\lambda \sum_{j=1}^{l_\lambda} (h_{(1)}^{\nu,i})_{(1)}^{\lambda,j}\otimes (h_{(1)}^{\nu,i})_{(2)}^{\lambda,j}\right)\otimes\sum_{i=1}^{n_\nu} h_{(2)}^{\nu,i}\right),
\end{multline}
здесь
$$
\Delta(h_{(1)}^{\nu,i})=\lim_\lambda \sum_{j=1}^{l_\lambda} (h_{(1)}^{\nu,i})_{(1)}^{\lambda,j}\otimes (h_{(1)}^{\nu,i})_{(2)}^{\lambda,j}\quad
\text{и}\quad
\Delta(h_{(2)}^{\nu,i})=\lim_\mu \sum_{j=1}^{m_\mu} (h_{(2)}^{\nu,i})_{(1)}^{\mu,j}\otimes (h_{(2)}^{\nu,i})_{(2)}^{\mu,j}.
$$

Согласно теореме о повторном пределе \cite[глава~2, теорема~4]{Kel} оба повторных предела в~\eqref{detcoass} могут быть заменены на простой. Это означает, в частности,  что итерацию коумножения можно обозначать через
$$
\sum h_{(1)}\otimes\cdots \otimes h_{(n)},
$$
где также подразумевается простой предел.

\subsection*{Смэш-произведения}
Пусть $H$ --- $\mathbin{\widehat{\otimes}}$-биалгебра и пусть $\mathbin{\widehat{\otimes}}$-алгебра~$A$  снабжена структурой левого $H$-$\mathbin{\widehat{\otimes}}$-модуля.
Тогда $A\mathbin{\widehat{\otimes}} A$ и $\mathbb{C}$ являются левыми $H$-$\mathbin{\widehat{\otimes}}$-модулями относительно внешних умножений, заданными формулами
$$
h\cdot(a\otimes b)\!:=\sum (h_{(1)}\cdot a)\otimes  (h_{(2)}\cdot b)
\qquad (h\in H,\,a,b\in A),
$$
$$
h\cdot\lambda\!:=\varepsilon(h)\lambda, \qquad (h\in H,\,\lambda\in\mathbb{C}).
$$

Напомним, что~$A$ называется (левой)  \emph{$H$-$\mathbin{\widehat{\otimes}}$-модульной алгеброй}, если морфизм умножения $\mu\!:A\mathbin{\widehat{\otimes}} A\to A$ и единичное отображение $\mathbb{C}\to A$ являются морфизмами левых $H$-модулей. Используя обобщенные обозначения Свидлера, запишем эти условия как
\begin{equation}\label{Hmodalgcond}
h\cdot(ab)=\sum (h_{(1)}\cdot a)(h_{(2)}\cdot b)\quad\text{и}\quad h\cdot 1=\varepsilon(h)1\qquad(h\in H,\,a,b\in A).
\end{equation}
(Чисто алгебраическую версию см. в \cite[\S\,11.2]{Ra}.)

Cмэш-произведения могут быть определены через явную конструкцию или через универсальное свойство. Начнём с универсального определения (в чисто алгебраическом случае оно введено в \cite{HS69}, см. также \cite{Mo77,Ma90a}).

Мы дополнительно предполагаем, что $H$ является $\mathbin{\widehat{\otimes}}$-алгеброй Хопфа. Пусть $\psi\!:H\to B$ --- гомоморфизм $\mathbin{\widehat{\otimes}}$-алгебр. Аналогично \cite{HS69}, зададим
\emph{присоединённое действие} $H$ на $B$ формулой
$$
h\cdot b \!:=\sum \psi(h_{(1)})\,b\,\psi(S(h_{(2)}))\qquad (h\in H,\,b\in B).
$$
Из \eqref{Hmodalgcond} и аксиомы антипода нетрудно вывести, что $B$ является $H$-$\mathbin{\widehat{\otimes}}$-модульной алгеброй относительно присоединённого действия.

\begin{df}\label{PsmaDef}
Пусть $H$ есть $\mathbin{\widehat{\otimes}}$-алгебра Хопфа, и $A$ является $H$-$\mathbin{\widehat{\otimes}}$-модульной алгеброй.  \emph{Аналитическим  смэш-произведением}  $A\mathop{\widehat{\#}} H$ называется  $\mathbin{\widehat{\otimes}}$-алгебра, снабжённая гомоморфизмами $\mathbin{\widehat{\otimes}}$-алгебр
$i\!:A\to A\mathop{\widehat{\#}} H$ и $j\!:  H\to A\mathop{\widehat{\#}} H$, такими что выполнены следующие свойства.

(A)~$i$  является гомоморфизмом $H$-$\mathbin{\widehat{\otimes}}$-модульных алгебр (относительно присоединёного действия, порождённого~$j$).

(B)~ Для любых $\mathbin{\widehat{\otimes}}$-алгебры $B$, гомоморфизма $\mathbin{\widehat{\otimes}}$-алгебр $\psi\!:  H\to  B$  и гомоморфизма $H$-$\mathbin{\widehat{\otimes}}$-модульных алгебр $\varphi\!: A\to B$ (относительно  присоединёного действия, порождённого~$\psi$) найдётся единственный гомоморфизм
$\mathbin{\widehat{\otimes}}$-алгебр $\theta$, такой что диаграмма
\begin{equation*}
   \xymatrix{
 & A\mathop{\widehat{\#}} H \ar@{-->}[dd]^{\theta}& \\
A \ar[ur]^i \ar[dr]_\varphi&&H\ar[ul]_j\ar[dl]^\psi\\
& B &}
\end{equation*}
коммутативна.
\end{df}

Явная конструкция (в аналитическом случае введенная в \cite{Pi4}) в обобщенных обозначениях Свидлера выглядит так же, как и в чисто алгебраическом случае. А именно, на полном локальном выпуклом пространстве   $A\mathbin{\widehat{\otimes}}H$ формула
\begin{equation}\label{exmusmpr}
(a\otimes h)(a'\otimes h')=\sum a(h_{(1)}\cdot a')\otimes h_{(2)}h'\qquad(h,h'\in H;\,a,a'\in A).
\end{equation}
задаёт непрерывное  ассоциативное умножение. Нетрудно видеть, что полученная алгебра и отображения $i\!:a\mapsto a\otimes 1$, $j\!:h\mapsto 1\otimes h$ удовлетворяют условиям из определения~\ref{PsmaDef}. Таким образом, подлежащее локально выпуклое пространство $\mathbin{\widehat{\otimes}}$-алгебры $A\mathop{\widehat{\#}}H$ совпадает с $A\mathbin{\widehat{\otimes}}H$.

\begin{rem}
Понятие полупрямого произведения групп может быть обобщено на случай алгебр Хопфа несколькими способами. Используемая здесь конструкция смэш-произведения, по-видимому, является простейшим из них. В аналитическом случае одно из более общих понятий рассмотрено в \cite[\S\,5]{Pi3} (в отличие от  \cite{Pi4} оно введено через универсальное свойство).
Его алгебраическая версия была известна ранее, например, она упоминается в \cite[\S\,8]{Ta80}, но приведённые там ссылки не проясняют происхождение понятия, так что можно считать, что оно относится к сфере математического фольклора.
\end{rem}

Нам также будут нужны достаточные условия того, что аналитическое смэш-произведение является $\mathbin{\widehat{\otimes}}$-алгеброй Хопфа.
Пусть $C$ --- $\mathbin{\widehat{\otimes}}$-коалгебра, снабжённая структурой левого $H$-$\mathbin{\widehat{\otimes}}$-модуля. Она называется (левой) \emph{$H$-$\mathbin{\widehat{\otimes}}$-модульной коалгеброй}, если коумножение $C\to C\mathbin{\widehat{\otimes}} СC$ и коединица $C\to \mathbb{C}$ являются морфизмами левых $H$-$\mathbin{\widehat{\otimes}}$-модулей (ср. \cite[Definition 2.1(c)]{Mo77}).  Если $A$ --- $\mathbin{\widehat{\otimes}}$-биалгебра, которая является
$H$-$\mathbin{\widehat{\otimes}}$-модульной алгеброй и $H$-$\mathbin{\widehat{\otimes}}$-модульной коалгеброй одновременно, то она
называется (левой)  \emph{$H$-$\mathbin{\widehat{\otimes}}$-модульной биалгеброй} (ср. \cite[Definition 2.1(e)]{Mo77}).

Тензорное произведение $\mathbin{\widehat{\otimes}}$-коалгебр $C$ и $C'$ определяется как локально выпуклое пространство $C\mathbin{\widehat{\otimes}} C'$ с коумножением
и коединицей, заданными формулами
\begin{equation}\label{tenptco}
c\otimes c' \mapsto \sum (c_{(1)}\otimes c'_{(1)}) \otimes  (c_{(2)}\otimes c'_{(2)})\quad\text{и}\quad c\otimes c' \mapsto \varepsilon_C(c)\varepsilon_{C'}(c').
\end{equation}

В следующем утверждении основным является дополнительное предположение о том, что $H$ кокоммутативна.

\begin{pr}\label{cocsmH}
Пусть $H$ --- кокоммутативная $\mathbin{\widehat{\otimes}}$-алгебра Хопфа, и $A$ есть $\mathbin{\widehat{\otimes}}$-алгебра Хопфа, которая является $H$-$\mathbin{\widehat{\otimes}}$-модульной биалгеброй. Рассмотрим $A\mathop{\widehat{\#}} H$ как тензорное произведение $\mathbin{\widehat{\otimes}}$-коалгебр.  Тогда коумножение совместимо со структурой $\mathbin{\widehat{\otimes}}$-алгебры и тем самым $A\mathop{\widehat{\#}} H$  является $\mathbin{\widehat{\otimes}}$-биалгеброй. Более того, она является $\mathbin{\widehat{\otimes}}$-алгеброй Хопфа относительно  антипода, заданного формулой
$$
a\otimes h \mapsto \sum S_H(h_{(2)})\cdot S_A(a)\otimes S_H(h_{(1)})
\qquad(h\in H,\,a\in A).
$$
\end{pr}
Доказательство в точности такое же, как в \cite[Theorem 2.13]{Mo77}) для ''обычных'' алгебр Хопфа и поэтому опущено.

Отметим, что все $\mathbin{\widehat{\otimes}}$-алгебры Хопфа, рассматриваемые в этой статье, кокоммутативны.

\subsection*{Смэш-произведения и алгебры Аренса-Майкла}
Также нам понадобится следующее понятие, которое может рассматриваться как частный случай предложенного Пирковским в \cite[Definition 5.5]{Pi3} определения смэш-произведения для алгебр Аренса-Майкла. Напомним, что \emph{алгеброй Аренса-Майкла} называется $\mathbin{\widehat{\otimes}}$-алгебра, топология на которой может быть задана семейством субмультипликативных преднорм, другими словами, она является проективным пределом системы банаховых алгебр.
Если $A$ --- ассоциативная топологическая алгебра, то её пополнение относительно семейства
всех субмультипликативных непрерывных преднорм называется её \emph{оболочкой Аренса-Майкла}. Она обозначается через $\widehat A$ и обладает следующим универсальным свойством: всякий непрерывный гомоморфизм в банахову алгебру может быть пропущен через $\widehat A$, ср. \S\,\ref{s:UNAL}. Подробности см. в \cite[глава 5]{X2}.

\begin{df}
Если в дополнение к условиям в определении~\ref{PsmaDef} $A$ и $H$ являются алгебрами Аренса-Майкла, то $A\mathop{\#}\nolimits_{\mathsf{AM}}H$ есть алгебра Аренса-Майкла, имеющая аналогичное универсальное свойство при условии, что $B$ --- также алгебра Аренса-Майкла.
\end{df}

Так же как и в \cite{Pi3} существование $A\mathop{\#}\nolimits_{\mathsf{AM}}H$ следует из универсального свойства оболочки Аренса-Майкла. А именно, мы можем положить
\begin{equation}\label{AMSmp}
A\mathop{\#}\nolimits_{\mathsf{AM}}H\!:=(A\mathop{\widehat{\#}} H)\sphat\,,
\end{equation}
где последний циркумфлекс обозначает оболочку Аренса-Майкла.

\section{Разложение алгебры аналитических функционалов полупрямого произведения}
\label{s:DAAF}

В этом разделе показано, что разложение комплексной группы Ли в полупрямое произведение порождает
разложение алгебры аналитических функционалов в аналитическое смэш-произведение.

Для комплексного многообразия $M$ мы обозначаем через $\mathcal{O}(M)$
локально выпуклое пространство  всех голоморфных функций на~$M$ (с топологией равномерной сходимости на компактных подмножествах), а через  ${\mathscr A}(M)$ --- сильное дуальное пространство $\mathcal{O}(M)'$  (множество непрерывных линейных функционалов с топологией равномерной сходимости на ограниченных подмножествах). Если $N$ и $M$ --- комплексные многообразия, то поскольку $\mathcal{O}(N)$ и
$\mathcal{O}(M)$ являются ядерными пространствами Фреше, имеют место следующие хорошо известные топологические изоморфизмы  полных локально выпуклых пространств (они понадобятся ниже):
\begin{multline}
\label{LSCtnpr}
 {\mathscr A}(N\times M)=\mathcal{O}(N\times M)'\cong\\
(\mathcal{O}(N)\mathbin{\widehat{\otimes}}\mathcal{O}(M))'\cong  \mathcal{O}(N)'\mathbin{\widehat{\otimes}}\mathcal{O}(M)'= {\mathscr A}(N)\mathbin{\widehat{\otimes}} {\mathscr A}(M).
\end{multline}

Пусть теперь $G$  --- комплексная группа Ли. Тогда $\mathcal{O}(G)$ является
$\mathbin{\widehat{\otimes}}$-алгеброй  (относительно поточечного умножения).
Формулы
\begin{equation*}
 \Delta_\mathcal{O} (f)(g, h) = f(gh),\quad \varepsilon_\mathcal{O} (f) = f(1),\quad  (S_\mathcal{O} f)(g) = f(g^{-1})
\end{equation*}
задают на ней каноническую структуру $\mathbin{\widehat{\otimes}}$-алгебры Хопфа.
Пространство ${\mathscr A}(G)$ можно снабдить дуальной структурой $\mathbin{\widehat{\otimes}}$-алгебры Хопфа. А именно, умножение (свёртка) на ${\mathscr A}(G)$  задаётся формулой
$$
\langle a'_1a'_2, f\rangle\!:= (a'_1\otimes a'_2)\Delta_\mathcal{O} (f) \qquad(a'_1,a'_2\in{\mathscr A}(G),\,f\in\mathcal{O}(G)),
$$
а единицей является дельта-функция в единице группы. Остальные операции определяются формулами
\begin{equation}\label{AGoper}
\Delta_{\mathscr A}(a')(f_1\otimes f_2)=\langle a', f_1f_2\rangle,\quad\varepsilon_{\mathscr A}(a')=\langle a', 1\rangle,\quad \langle S_{\mathscr A}(a'), f\rangle = \langle a', S_\mathcal{O}(f)\rangle
\end{equation}
(здесь $a'\in{\mathscr A}(G)$ и $f,f_1,f_2\in\mathcal{O}(G)$).
Очевидно, что ${\mathscr A}(G)$  кокоммутативна.

Пусть $G_1\rtimes G_2$ --- полупрямое произведение комплексных групп Ли
относительно гомоморфизма $\alpha\!:G_2\to{\mathop{\mathrm{Aut}}\nolimits}(G_1)$. Рассмотрим кодействие
$$
\beta\!:\mathcal{O}(G_1)\to \mathcal{O}(G_2\times G_1);\quad (\beta
f)(g_2,g_1)\!:=f(\alpha_{g_2}(g_1))\,.
$$
и (отождествляя $\mathcal{O}(G_2\times G_1)$ с $\mathcal{O}(G_2)\mathbin{\widehat{\otimes}} \mathcal{O}(G_1))$ дуальное действие ${\mathscr A}(G_2)$ на ${\mathscr A}(G_1)$:
\begin{equation}\label{actAG2}
(b'\cdot a')(f)\!:=(b'\otimes a')(\beta( f))\qquad (b'\in
{\mathscr A}(G_2),\,a'\in {\mathscr A}(G_1),\,f\in \mathcal{O}(G_1)).
\end{equation}

\begin{thm} \label{AGsemid}
Пусть $G_1\rtimes G_2$ --- полупрямое произведение комплексных групп Ли
относительно гомоморфизма $\alpha\!:G_2\to{\mathop{\mathrm{Aut}}\nolimits}(G_1)$. Тогда ${\mathscr A}(G_1)$ является
левой ${\mathscr A}(G_2)$-$\mathbin{\widehat{\otimes}}$-модульной
биалгеброй относительно действия, заданного~\eqref{actAG2}, а топологический изоморфизм локально выпуклых пространств
$$
\varphi\!:{\mathscr A}(G_1)\mathop{\widehat{\#}} {\mathscr A}(G_2)\to{\mathscr A}(G_1\rtimes G_2)
$$
\emph{(}см.~\eqref{LSCtnpr}\emph{)} является изоморфизмом
 $\mathbin{\widehat{\otimes}}$-алгебр Хопфа.
\end{thm}

\begin{proof}
Рассмотрим очевидно определённые гомоморфизмы из групповых алгебр $\mathbb{C}[G_1]$ и $\mathbb{C}[G_2]$ в соответствующие алгебры аналитических функционалов. Заметим, что их образы плотны. (Действительно, в силу теоремы Хана-Банаха достаточно убедиться, что всякий непрерывный линейный функционал на алгебре аналитических функционалов, значения которого равны~$0$  на всех дельта-функциях, тождественно равен~$0$. Последнее очевидно, поскольку непрерывный линейный функционал определяется голоморфной функцией на группе ввиду рефлексивности пространства голоморфных функций.) Так как $\mathbb{C}[G_1]$ является $\mathbb{C}[G_2]$-модульной биалгеброй, то в силу плотности образов указанных гомоморфизмов аксиомы $\mathbin{\widehat{\otimes}}$-модульной биалгебры выполнены для ${\mathscr A}(G_1)$, поскольку они выполнены для $\mathbb{C}[G_1]$. Таким образом, ${\mathscr A}(G_1)\mathop{\widehat{\#}} {\mathscr A}(G_2)$ корректно определено и является $\mathbin{\widehat{\otimes}}$-алгеброй Хопфа, так как
${\mathscr A}(G_2)$ кокоммутативна, см. предложение~\ref{cocsmH}.

Чтобы убедиться в том, что отображение $\varphi$ есть гомоморфизм $\mathbin{\widehat{\otimes}}$-алгебр Хопфа достаточно установить, что оно является гомоморфизмом $\mathbin{\widehat{\otimes}}$-алгебр и  гомоморфизмом $\mathbin{\widehat{\otimes}}$-коалгебр. Из плотности образов $\mathbb{C}[G_1]$ и $\mathbb{C}[G_2]$ следует, что образ $\mathbb{C}[G_1]\mathop{\#} \mathbb{C}[G_2]$ плотен в ${\mathscr A}(G_1)\mathop{\widehat{\#}} {\mathscr A}(G_2)$. Поэтому совместимость $\varphi$ с соответствующими операциями достаточно проверить только для элементов вида $\delta_{g_1}\otimes \delta_{g_2}$. Заметим, что $\varphi(\delta_{g_1}\otimes \delta_{g_2})= \delta_{(g_1,g_2)}$. Тогда из формулы умножения~\eqref{exmusmpr} и формулы умножения в полупрямом произведении групп следует, что $\varphi$  мультипликативно.  Из формулы коумножения  \eqref{tenptco} следует, что $(\varphi\otimes \varphi)(\Delta_{\#}(\delta_{g_1}\otimes \delta_{g_2}))= \Delta (\varphi(\delta_{g_1}\otimes \delta_{g_2}))$ (здесь $\Delta_{\#}$ --- коумножение в ${\mathscr A}(G_1)\mathop{\widehat{\#}} {\mathscr A}(G_2)$, а $\Delta$ --- в ${\mathscr A}(G_1\rtimes G_2)$). Отсюда получаем, что  $\varphi$ комультипликативно. Для единицы и коединицы утверждения тривиальны.

Будучи топологическим изоморфизмом локально выпуклых пространств, $\varphi$ тем самым является изоморфизмом $\mathbin{\widehat{\otimes}}$-алгебр Хопфа.
\end{proof}

\begin{co}\label{GDAF}
Если комплексная группа Ли $G$ представлена как итерированное полупрямое произведение:
$$
G=((\cdots (F_1 \rtimes F_2)\rtimes\cdots)
\rtimes F_n)
$$
для некоторых комплексных групп Ли $F_1,\ldots, F_n$, то алгебра аналитических функционалов
топологически изоморфна итерированному аналитическому смэш-произ\-ве\-дению:
\begin{equation}\label{itsmpr}
{\mathscr A}(G)\cong ((\cdots ({\mathscr A}(F_1) \mathop{\widehat{\#}} {\mathscr A}(F_2))\mathop{\widehat{\#}}\cdots)
\mathop{\widehat{\#}} {\mathscr A}(F_n)).
\end{equation}
\end{co}
\begin{rem}
Если $G$ линейна (существует инъективное конечномерное голоморфное представление) и связна, то её можно представить в виде $B\rtimes L$, где $B$ разрешима и односвязна, а $L$ комплексно линейно редуктивна (является комплексификацией компактной действительной группы Ли), см., например, \cite[с.\,601, Theorem 16.3.7]{HiNe}.
Тогда $B$  разлагается в итерированное полупрямое произведение нескольких экземпляров группы~$\mathbb{C}$
(ср. \cite[с.\,449, Theorem 11.2.14]{HiNe}). Таким образом, все множители в \eqref{itsmpr}, кроме последнего, будут изоморфны ${\mathscr A}(\mathbb{C})$, преобразование Фурье которой, как известно, есть алгебра целых функций экспоненицального типа на~$\mathbb{C}$, см. например, \cite[предложение 7.2]{Ak08}, а последний множитель ${\mathscr A}(L)$ полностью определяется теорией представлений~$L$.
\end{rem}

\section{Разложение пополнений, ассоциированных с субмультипликативными весами}
\label{s:DCASW}

В этом разделе найдены достаточные условия того, что разложение комплексной группы Ли в полупрямое произведение порождает
разложение пополнения алгебры аналитических функционалов, ассоциированного с субмультипликативным весом, в аналитическое смэш-произведение.

\subsection*{Пополнения пространств аналитических функционалов}
Следуя \cite{Ak08}, обозначим для комплексно-аналитического многообразия $M$ и локально ограниченной функции $\upsilon\!:M\to [1,+\infty)$ замыкание в~${\mathscr A}(M)$ абсолютно выпуклой оболочки подмножества
\begin{equation}\label{Vup}
\{\upsilon(x)^{-1}\delta_x:\,x\in M\}
\end{equation}
через $V_\upsilon$ (здесь $\delta_x$ обозначает дельта-функцию в точке $x$).  Функционал Минковского на ${{\mathscr A}}(M)$, ассоциированный с $V_\upsilon$, является непрерывной преднормой на ${{\mathscr A}}(M)$. Обозначим его через $\|\cdot\|_{\upsilon}$,  через  ${{\mathscr A}}_{\upsilon}(M)$ --- пополнение ${{\mathscr A}}(M)$ относительно $\|\cdot\|_{\upsilon}$, а через
${{\mathscr A}}_{\upsilon^\infty}(M)$ --- пополнение ${{\mathscr A}}(M)$ относительно последовательности преднорм $$
(\|\cdot\|_{\upsilon^n};\,n\in\mathbb{N}),\qquad \text{где $\upsilon^n(x)\!:=\upsilon(x)^n$.}
$$

\begin{lm}\label{tenho}
Пусть $\upsilon_1$ и $\upsilon_2$ --- локально ограниченные функции на
комплексно-ана\-ли\-ти\-ческом многообразии $M$ и пусть $\upsilon(x)\!:=\upsilon_1(x)\upsilon_2(x)$. Тогда диагональное вложение $M\to M\times  M$  индуцирует непрерывное линейное отображение
$$
{\mathscr A}_{\upsilon}(M)\to{\mathscr A}_{\upsilon_1}(M)\mathbin{\widehat{\otimes}} {\mathscr A}_{\upsilon_2}(M).
$$
\end{lm}
\begin{proof}
Рассмотрим на $M\times M$ функцию
$$
\bar\upsilon(x,y)\!:=\upsilon_1(x)\upsilon_2(y).
$$
Согласно  \cite[Proposition~5.5]{ArAnF}, банахово пространство ${\mathscr A}_{\upsilon_1}(M)\mathbin{\widehat{\otimes}} {\mathscr A}_{\upsilon_2}(M)$
изометрически изоморфно ${\mathscr A}_{\bar\upsilon}(M\times M)$.
В силу определения, единичные шары $V_{\upsilon}$   и $V_{\bar\upsilon}$ в ${\mathscr A}(M)$ и ${\mathscr A}(M\times M)$  есть замыкания абсолютно выпуклых оболочек  множеств
\begin{equation*}
V_{\upsilon}^0\!:=\{\upsilon(x)^{-1}\delta_x\!:\,x\in M\}\qquad\text{и}\qquad
V_{\bar\upsilon}^0\!:=\{\bar{\upsilon}(x,y)^{-1}\delta_{x}\otimes \delta_{y}\!:\,x,y\in M\}
\end{equation*}
относительно соответственно $\|\cdot\|_{\upsilon}$ и $\|\cdot\|_{\bar\upsilon}$.

Обозначим через $\Delta$ отображение ${\mathscr A}(M)\to {\mathscr A}(M\times M)$, индуцированное диагональным вложением (с точностью до изоморфизма ${\mathscr A}(M)\mathbin{\widehat{\otimes}} {\mathscr A}(M)\cong{\mathscr A}(M\times M)$, см. \eqref{LSCtnpr}).
Пусть $a'$ есть конечная сумма $\sum c_x\,\upsilon(x)^{-1}\delta_{x}$, где $\sum |c_x|\le 1$.
Так как $\Delta(\delta_x)=\delta_x\otimes\delta_x$ для каждого $x\in M$, то
$$
\Delta(a')=\sum_x c_x\,\upsilon(x)^{-1}\,\delta_x\otimes\delta_x=\sum_{x,y} c_{xy}\,\bar\upsilon(x,y)^{-1}\,\delta_x\otimes \delta_y,
$$
где $c_{xx}\!:=c_{x}$ и $c_{xy}\!:=0$, если $x\ne y$. Таким образом, $\Delta$ отображает абсолютно выпуклую оболочку множества
$V_{\upsilon}^0$  в абсолютно выпуклую оболочку $V_{\bar\upsilon}^0$.

Переходя к замыканиям, заметим, что поскольку $\Delta$  непрерывно, образ множества $V_{{\upsilon}}$ содержится в $V_{\bar\upsilon}$. Таким образом, ${\mathscr A}_{\upsilon}(M)\to{\mathscr A}_{\bar\upsilon}(M\times M)$ является сжимающим, в частности, оно непрерывно.
\end{proof}

\subsection*{Cубмультипликативные веса и полупрямые произведения}
Строго положительная локально ограниченная
функция $\omega\!: G \to\mathbb{R}$  на локально компактной группе $G$ называется \emph{субмультипликативным весом}\footnote{Описанный здесь подход близок к первоначальным определениям из \cite{Ak08}. В работах автора \cite{ArAMN,ArAnF,ArHR,AHHFG,ArUAPI} вместо весов часто используется эквивалентное описание через \emph{функции длины}  (задаются условием субаддитивности). Переход от субмультипликативного веса к функции длины осуществляется логарифмированием.}, если $\omega(gh)\le \omega(g)\omega(h)$ для всех $g, h \in G$.  Нас интересуют субмультипликативные веса на комплексной группе Ли $G$ и пополнения алгебры ${\mathscr A}(G)$  вида ${\mathscr A}_{\omega^\infty}(G)$.

\begin{pr}\label{ominftyAM}
Пусть $\omega$ --- субмультипликативный вес на  комплексной группе Ли $G$, такой что $\omega(g)\ge 1$ для всех $g\in G$. Тогда умножение на ${\mathscr A}(G)$ продолжается до умножения на ${\mathscr A}_{\omega^\infty}(G)$, превращающего её в алгебру Фреше-Аренса-Майкла. Более того, ${\mathscr A}_{\omega^\infty}(G)$ является кокоммутативной $\mathbin{\widehat{\otimes}}$-биалгеброй относительно операций, непрерывно продолженных с ${\mathscr A}(G)$, а отображение пополнения ${\mathscr A}(G)\to {\mathscr A}_{\omega^\infty}(G)$ является гомоморфизмом $\mathbin{\widehat{\otimes}}$-биалгебр.
\end{pr}
\begin{proof}
Так как топология на пространстве ${\mathscr A}_{\omega^\infty}(G)$ задана
счётным семейством преднорм $(\|\cdot\|_{\omega^n})$, оно является пространством Фреше. Заметим, что
$\omega^n$ --- субмультипликативный вес для каждого~$n$ и тем самым
$\|\cdot\|_{\omega^n}$ субмультипликативна и непрерывна на
${\mathscr A}(G)$ согласно \cite[лемма 5.1(1)]{Ak08}. Следовательно, ${\mathscr A}_{\omega^\infty}(G)$ является алгеброй Аренса-Майкла,
а ${\mathscr A}(G)\to {\mathscr A}_{\omega^\infty}(G)$ --- непрерывный гомоморфизм $\mathbin{\widehat{\otimes}}$-алгебр.

В силу леммы~\ref{tenho} диагональное вложение $G\to G\times  G$ (и соответствующее ему коумножение $\Delta_{\mathscr A}\!:{\mathscr A}(G)\to {\mathscr A}(G\times G)$) для каждого $n\in\mathbb{N}$ индуцирует непрерывное линейное отображение
${\mathscr A}_{\omega^{2n}}(G)\to {\mathscr A}_{\omega^n}(G)\mathbin{\widehat{\otimes}} {\mathscr A}_{\omega^n}(G)$. В пределе, когда $n\to \infty$,  мы получаем отображение ${\mathscr A}_{\omega^\infty}(G)\to {\mathscr A}_{\omega^\infty}(G)\mathbin{\widehat{\otimes}} {\mathscr A}_{\omega^\infty}(G)$.

Чтобы убедиться, что коединица продолжается на ${\mathscr A}_{\omega^\infty}(G)$, достаточно проверить, что $|\varepsilon(a')|\le \|a'\|_{\omega}$, где $a'$ --- конечная сумма $\sum c_g\,\omega(g)^{-1}\delta_{g}$, такая что $\sum |c_g|\le 1$. Из определения $\varepsilon$ (см.~\eqref{AGoper}) следует, что $\varepsilon(a')=\sum c_g\,\omega(g)^{-1}$.
Поскольку $\omega(g)\ge 1$, имеем $|\varepsilon(a')|\le \sum |c_g|\le 1$. Итак, мы  получили непрерывное отображение, индуцированное коединицей.

Аксиомы кокоммутативной $\mathbin{\widehat{\otimes}}$-биалгебры для ${\mathscr A}_{\omega^\infty}(G)$ выполнены в силу того, что они выполнены для ${\mathscr A}(G)$, непрерывности операций и плотности ${\mathscr A}(G)$ в ${\mathscr A}_{\omega^\infty}(G)$ (ср. \cite[Lemma 5.1]{AHHFG}). При этом автоматически получаем, что ${\mathscr A}(G)\to {\mathscr A}_{\omega^\infty}(G)$ является гомоморфизмом $\mathbin{\widehat{\otimes}}$-биалгебр.
\end{proof}

\begin{rem}\label{lowestw}
В предложении~\ref{ominftyAM} условие на значения веса может быть заменено на формально более слабое: найдётся $c>0$, такое что $\omega\ge c$. Доказательство при этом не изменится. Однако это излишне, поскольку, как нетрудно показать, указанное условие влечёт, что $\omega\ge 1$, см., например, \cite[Proposition 2.2]{Dz86}.
\end{rem}

Пусть $\omega_1$ и $\omega_2$ --- субмультипликативные веса на $G$. Мы говорим, что $\omega_2$ \emph{мажорирует} $\omega_1$ (пишем
$\omega_1\preceq \omega_2$), если найдутся $C>0$ и $\gamma>0$, такие что
$$
\omega_1(g)\le C\,\omega_2(g)^\gamma \quad\text{для всех $g\in G$.}
$$
Если $\omega_1\preceq \omega_2$ и $\omega_2\preceq \omega_1$, то говорим, что $\omega_1$ и $\omega_2$ \emph{эквивалентны} (пишем
$\omega_1\simeq\omega_2$).

Нетрудно видеть, что ${\mathscr A}_\omega(G)$ является банаховой алгеброй для всякого субмультипликативного веса на комплексной группе Ли~$G$. Эквивалентные веса порождают изоморфные алгебры Фреше вида ${\mathscr A}_{\omega^\infty}(G)$, однако, как видно из следующего примера,  могут порождать неизоморфные банаховы алгебры вида ${\mathscr A}_\omega(G)$.
\begin{exm}
Веса $\omega(z)=1+|z|$ и $\omega_1(z)=\sqrt{1+|z|}$ на $\mathbb{C}$ эквивалентны. Поэтому $\mathbin{\widehat{\otimes}}$-алгебры Хопфа ${\mathscr A}_{\omega^\infty}(\mathbb{C})$ и  ${\mathscr A}_{\omega_1^\infty}(\mathbb{C})$ изоморфны. Однако банаховы алгебры ${\mathscr A}_{\omega}(\mathbb{C})$ и  ${\mathscr A}_{\omega_1}(\mathbb{C})$ не являются изоморфными. Это следует из того, что двойственные банаховы пространства $\mathcal{O}_\omega(\mathbb{C})$ и $\mathcal{O}_{\omega_1}(\mathbb{C})$  (см. определение в~\eqref{cOom} ниже) не изоморфны.  А~именно, $\mathcal{O}_\omega(\mathbb{C})$ состоит из  многочленов степени не выше~$1$, а $\mathcal{O}_{\omega_1}(\mathbb{C})$ --- из констант. Отсюда нетрудно вывести, что ${\mathscr A}_\omega(\mathbb{C})\cong\mathbb{C}[x]/(x^2)$, а  ${\mathscr A}_{\omega_1}(\mathbb{C})\cong \mathbb{C}$, ср. доказательство леммы~\ref{1dimdsp}.
\end{exm}

Напомним, что субмультипликативный вес на локально компактной группе $G$
называется \emph{симметричным}, если $\omega(g^{-1})=\omega(g)$ для всех $g \in G$ и $\omega(1)=1$. Для наших целей пригоден и следующий ослабленный вариант этого понятия.
\begin{df}
Назовём субмультипликативный вес на локально компактной группе $G$ \emph{асимптотически симметричным}, если $\omega$ мажорирует функцию $g\mapsto \omega(g^{-1})$ и $\omega(g)\ge 1$ для всех $g \in G$.
\end{df}
Легко видеть, что всякий симметричный вес асимптотически симметричен.

\begin{pr}\label{ominftyAM2}
Пусть $\omega$ --- асимптотически симметричный субмультипликативный вес на  комплексной группе Ли $G$. Тогда ${\mathscr A}_{\omega^\infty}(G)$ является кокоммутативной $\mathbin{\widehat{\otimes}}$-алгеброй Хопфа относительно антипода, непрерывно продолженного с ${\mathscr A}(G)$, а отображение пополнения
${\mathscr A}(G)\to {\mathscr A}_{\omega^\infty}(G)$ является гомоморфизмом $\mathbin{\widehat{\otimes}}$-алгебр Хопфа.
\end{pr}
\begin{proof}
Согласно предложению~\ref{ominftyAM} ${\mathscr A}_{\omega^\infty}(G)$ является кокоммутативной $\mathbin{\widehat{\otimes}}$-биалгеброй. Осталось проверить утверждения, относящиеся к антиподу.

Так как $\omega$ асимптотически симметричен, найдутся $C>0$ и $\gamma>0$, такие что $\omega(g^{-1})\le C\,\omega(g)^\gamma$ для всех $g\in G$. В силу условия $\omega(g)\ge 1$, увеличивая $\gamma$ при необходимости, можно полагать, что $\gamma\in \mathbb{N}$.

Пусть $n \in \mathbb{N}$ и $a'$ есть конечная сумма $\sum c_g\,\omega(g)^{-n}\delta_{g}$, где $\sum |c_g|\le 1$. Так как $S(\delta_{g})=\delta_{g^{-1}}$, получаем
$$
S(a')=\sum c_{g}\,\omega(g)^{-n}\delta_{g^{-1}}=\sum c'_{g}\,\omega(g^{-1})^{-n/\gamma}\delta_{g^{-1}}=\sum c'_{g^{-1}}\,\omega(g)^{-n/\gamma}\delta_{g},
$$
где $c'_{g}\!:=c_{g}\,\omega(g^{-1})^{n/\gamma}\omega(g)^{-n}$. Так как
$$
\sum |c'_{g^{-1}}|\le \max_g\left\{|\omega(g)^{n/\gamma}\omega(g^{-1})^{-n}|\right\}\,\sum |c_{g^{-1}}|\le C^{n/\gamma},
$$
то $ S(a')\in C^{n/\gamma}V_{\omega^{n/\gamma}}$ (см.~\eqref{Vup}). Таким образом, $S$ порождает линейное непрерывное отображение ${\mathscr A}_{\omega^n}(G)\to {\mathscr A}_{\omega^{n/\gamma}}(G)$. Переходя к пределу по $n$, кратным $\gamma$, получаем линейное непрерывное отображение ${\mathscr A}_{\omega^\infty}(G)\to {\mathscr A}_{\omega^\infty}(G)$. Аксиомы антипода выполнены в силу плотности ${\mathscr A}(G)$ в ${\mathscr A}_{\omega^\infty}(G)$

Всякий гомоморфизм $\mathbin{\widehat{\otimes}}$-биалгебр между $\mathbin{\widehat{\otimes}}$-алгебрами Хопфа является гомоморфизмом $\mathbin{\widehat{\otimes}}$-алгебр Хопфа (доказательство этого факта такое же, как в алгебраческом случае, см. \cite[p.\,152, Proposition 4.2.5]{DNR}). Значит, в силу предложения~\ref{ominftyAM} это верно и для
${\mathscr A}(G)\to {\mathscr A}_{\omega^\infty}(G)$.
\end{proof}

\begin{pr}\label{moalgwie}
Пусть $G_1$ и $G_2$ --- комплексные группы Ли, a гомоморфизм $\alpha\!:G_2\to{\mathop{\mathrm{Aut}}\nolimits}(G_1)$ задаёт голоморфное действие на~$G_1$.  Если субмультипликативные веса $\omega_1$ и $\omega_2$ на соответственно $G_1$ и $G_2$ таковы, что
$$
\omega_1(\alpha_{g_2}(g_1))\preceq \omega_1(g_1)\omega_2(g_2)\quad\text{на
$\,G_1\times G_2$} \quad(\text{$g_1\in G_1$, $g_2\in G_2$}),
$$
то ${\mathscr A}_{\omega_1^\infty}(G_1)$ является  ${\mathscr A}_{\omega_2^\infty}(G_2)$-$\mathbin{\widehat{\otimes}}$-модульной биалгеброй.
\end{pr}
\begin{proof}
Достаточно показать, что действие ${\mathscr A}(G_2)\times {\mathscr A}(G_1) \to {\mathscr A}(G_1)$, индуцированное отображением $\alpha$, непрерывно в топологиях, порожденных последовательностями преднорм $(\|\cdot\|_{\omega_i^n}:, n\in \mathbb{N})$ на алгебрах ${\mathscr A}(G_i)$ ($i = 1, 2$). (Всё остальное следует из теоремы~\ref{AGsemid} и плотности образов ${\mathscr A}(G_1)$ и ${\mathscr A}(G_2)$ соответственно в  ${\mathscr A}_{\omega_1^\infty}(G_1)$ и ${\mathscr A}_{\omega_2^\infty}(G_2)$.)
Таким образом, нужно проверить, что для всякого $n\in\mathbb{N}$
найдутся $n_1,n_2\in\mathbb{N}$ и $K>0$, такие что
$$
\|b'\cdot a'\|_{\omega_1^n}\le K\|b'\|_{\omega_2^{n_2}} \,\|a'\|_{\omega_1^{n_1}}\quad\text{для всех $b'\in{\mathscr A}(G_2)$, $a'\in {\mathscr A}(G_1)$}.
$$

Пусть $C>0$ и $\gamma>0$ таковы, что
\begin{equation}\label{om12Cga}
\omega_1(\alpha_{g_2}(g_1))\le C(\omega_1(g_1)\,\omega_2(g_2))^\gamma \quad\text{для всех $g_1\in G_1$, $g_2\in G_2$.}
\end{equation}
Зафиксируем $n\in \mathbb{N}$ и рассмотрим на $G_1$ и $G_2$ аналитические функционалы
$$
a'=\sum_{g_1\in G_1} c'_{g_1}\omega_1(g_1)^{-n\gamma}\delta_{g_1}\quad\text{и} \quad
b'=\sum_{g_2\in G_2}  c''_{g_2}\omega_2(g_2)^{-n\gamma}\delta_{g_2},
$$
где  $\sum |c'_{g_1}|\le 1$, $\sum |c''_{g_2}|\le 1$, и множества коэффициентов, отличных от нуля, конечны. Тогда
$$
b'\cdot a'=\sum_{g_1,g_2} c'_{g_1}c''_{g_2} \omega_1(g_1)^{-n\gamma}\omega_2(g_2)^{-n\gamma} \delta_{\alpha_{g_2}(g_1)}=\sum_{h\in G_1} c_h \delta_h,
$$
где
$$
c_h=\sum_{h=\alpha_{g_2}(g_1)} c'_{g_1}c''_{g_2} \omega_1(g_1)^{-n\gamma}\omega_2(g_2)^{-n\gamma}.
$$
Положим $\widetilde c_h\!:=c_h \omega_1(h)^n$.
Из \eqref{om12Cga} имеем
$\omega_1(g_1)^{-n\gamma}\omega_2(g_2)^{-n\gamma}\omega_1(\alpha_{g_2}(g_1))^n  \le C^n$, поэтому
$$
|\widetilde c_h|\le  C^n \sum_{h=\alpha_{g_2}(g_1)} |c'_{g_1}|\,|c''_{g_2}|.
$$
Следовательно, $\sum_h |\widetilde c_h|\le C^n$.
Тем самым $b'\cdot a'\in C^n\,V_{\omega_1^n}$ (см.~\eqref{Vup}). Так как действие ${\mathscr A}(G_2)\times {\mathscr A}(G_1)\to {\mathscr A}(G_1)$, определённое в \eqref{actAG2}, непрерывно, отсюда следует, что оно отображает произведение единичных шаров $V_{\omega_2^{n\gamma}}\times V_{\omega_1^{n\gamma}}$ в шар $C^n\,V_{\omega_1^n}$. Итак, мы можем положить $n_1=n_2=n\gamma$ и $K=C^n$, что завершает рассуждение.
\end{proof}

\begin{thm}\label{smprde}
Предположим, что комплексная группа Ли $G_2$  голоморфно действует автоморфизмами на
комплексной группе Ли $G_1$.
Пусть $\omega_1$, $\omega_2$ и $\omega$  --- субмультипликативные веса на $G_1$, $G_2$ и $G_1\rtimes G_2$ соответственно, причём $\omega_1$ и $\omega_2$ асимптотически симметричны и
\begin{equation}\label{eqom1om2}
\omega(g_1g_2)\simeq  \omega_1(g_1)\omega_2(g_2)\quad\text{на
$\,G_1\times G_2$} \quad(\text{$g_1\in G_1$, $g_2\in G_2$}).
\end{equation}
Тогда выполнены следующие утверждения.

\emph{(A)}~${\mathscr A}_{\omega_1^\infty}(G_1)$ является  ${\mathscr A}_{\omega_2^\infty}(G_2)$-$\mathbin{\widehat{\otimes}}$-модульной биалгеброй.

\emph{(B)}~$\omega$ асимптотически симметричен.

\emph{(C)}~Изоморфизм из теоремы~\ref{AGsemid} индуцирует
изоморфизм $\mathbin{\widehat{\otimes}}$-алгебр Хопфа
$$
{\mathscr A}_{\omega_1^\infty}(G_1)\mathop{\widehat{\#}} {\mathscr A}_{\omega_2^\infty}(G_2)\to
{\mathscr A}_{\omega^\infty}(G_1\rtimes G_2);
$$
более того,
$$
{\mathscr A}_{\omega_1^\infty}(G_1)\mathop{\#}\nolimits_{\mathsf{AM}}  {\mathscr A}_{\omega_2^\infty}(G_2)\cong {\mathscr A}_{\omega_1^\infty}(G_1)\mathop{\widehat{\#}} {\mathscr A}_{\omega_2^\infty}(G_2).
$$
\end{thm}

\begin{proof}
(A)~Обозначим $G_1\rtimes G_2$ через $G$. Так как $\omega_2$ асимптотически симметричен, имеем
$$
\omega_1(g_2^{-1}g_1g_2)\simeq \omega(g_2^{-1}g_1g_2)\le \omega(g_2^{-1})\omega(g_1)\omega(g_2)\simeq \omega_2(g_2)^2\omega_1(g_1),
$$
где  $g_1\in G_1$ и $g_2\in G_2$ рассмотрены также как элементы $G$.
Кроме того, так как $\omega_1$ также асимптотически симметричен, то $\omega_1\ge 1$. Поэтому $\omega_2(g_2)^2\omega_1(g_1)\le \omega_2(g_2)^2\omega_1(g_1)^2$.
Следовательно,
\begin{equation}\label{om1al}
\omega_1(\alpha_{g_2}(g_1))\preceq \omega_1(g_1)\,\omega_2(g_2)\qquad\text{на
$G_1\times G_2$}.
\end{equation}
Применяя предложение~\ref{moalgwie}, получаем, что ${\mathscr A}_{\omega_1^\infty}(G_1)$ является ${\mathscr A}_{\omega_2^\infty}(G_2)$-$\mathbin{\widehat{\otimes}}$-модульной биалгеброй.

(B)~Пусть $g=g_1g_2$, где $g_1\in G_1$ и $g_2\in G_2$. Тогда из \eqref{eqom1om2} и \eqref{om1al} следует, что
$$
\omega(g^{-1})=\omega(\alpha_{g_2^{-1}}(g_1^{-1})g_2^{-1})\preceq \omega_1(\alpha_{g_2^{-1}}(g_1^{-1}))\,\omega_2(g_2^{-1})\preceq \omega_1(g_1^{-1})\,\omega_2(g_2^{-1}).
$$
Так как $\omega_1$ и $\omega_2$ асимптотически симметричны, то $\omega_1(g_1^{-1})\preceq\omega_1(g_1)$ и $\omega_2(g_2^{-1})\preceq\omega_2(g_2)$. Применяя \eqref{eqom1om2} в обратную сторону, получаем отсюда, что $\omega(g^{-1})\preceq\omega(g)$. С другой стороны, так как $\omega_1\ge 1$ и $\omega_2\ge 1$, то из \eqref{eqom1om2} следует, что $\{\omega(g):\,g\in G\}$ отделено от $0$  положительным числом. В силу \cite[Proposition 2.2]{Dz86} получаем, что $\omega(g)\ge 1$ для всех~$g$ (ср. замечание~\ref{lowestw}). Итак, $\omega$ асимптотически симметричен.

(C)~Так как ${\mathscr A}_{\omega_1^\infty}(G_1)$ является ${\mathscr A}_{\omega_2^\infty}(G_2)$-$\mathbin{\widehat{\otimes}}$-модульной биалгеброй, то корректно определено смэш-произведение ${\mathscr A}_{\omega_1^\infty}(G_1)\mathop{\widehat{\#}} {\mathscr A}_{\omega_2^\infty}(G_2)$, которое является $\mathbin{\widehat{\otimes}}$-алгеброй Хопфа  в силу предложений~\ref{ominftyAM2} и~\ref{cocsmH}. Поскольку $\omega$ асимптотически симметричен, из предложения~\ref{ominftyAM2}
следует, что ${\mathscr A}_{\omega^\infty}(G_1\rtimes G_2)$ также является $\mathbin{\widehat{\otimes}}$-алгеброй Хопфа.

Согласно \cite[Proposition~5.5]{ArAnF} из~\eqref{eqom1om2} следует, что топологический изоморфизм ${\mathscr A}(G_1)\mathop{\widehat{\#}}{\mathscr A}(G_2)\to{\mathscr A}(G_1\rtimes G_2)$ из теоремы~\ref{AGsemid} индуцирует топологический изоморфизм локально выпуклых пространств
$$
{\mathscr A}_{\omega_1^\infty}(G_1)\mathop{\widehat{\#}} {\mathscr A}_{\omega_2^\infty}(G_2)\to
{\mathscr A}_{\omega^\infty}(G_1\rtimes G_2).
$$
Так же, как в доказательстве теоремы~\ref{AGsemid}, из плотности образов групповых алгебр получаем, что он также является гомоморфизмом  $\mathbin{\widehat{\otimes}}$-алгебр Хопфа, а значит и их  изоморфизмом.

Далее, из~\eqref{AMSmp} следует, что ${\mathscr A}_{\omega_1^\infty}(G_1)\mathop{\#}\nolimits_{\mathsf{AM}}{\mathscr A}_{\omega_2^\infty}(G_2)$ изоморфна оболочке Аренса-Майкла алгебры ${\mathscr A}_{\omega_1^\infty}(G_1)\mathop{\widehat{\#}} {\mathscr A}_{\omega_2^\infty}(G_2)$, однако последняя, как только что показано, изоморфна ${\mathscr A}_{\omega^\infty}(G_1\rtimes G_2)$, которая есть алгебра Аренса-Майкла (предложение~\ref{ominftyAM}) и поэтому совпадает со своей оболочкой.
\end{proof}

\begin{exm}
Разложение вида~\eqref{eqom1om2} не всегда возможно. Рассмотрим на $\mathbb{C}\times \mathbb{C}$ субмультипликативный вес
$\omega\!:(u,v)\mapsto e^{|u+v|}$. Если $\omega(u,v)\simeq  \omega_1(u)\omega_2(v)$ для некоторых субмультипликативных весов $\omega_1$ и $\omega_2$, то, полагая $v=-u$, получаем, что функция $u\mapsto\omega_1(u)\omega_2(-u)$ ограничена. Если, кроме того $\omega_2$ симметричен, то функция $u\mapsto\omega_1(u)\omega_2(u)$ ограничена и тем самым $\omega(u,u)=e^{2|u|}$ ограничена. Противоречие.

Достаточные условия для разложения вида \eqref{eqom1om2} см., например, в \cite[Lemma~4.5]{ArAnF}.
\end{exm}

Нам понадобится следующая лемма, которая выводится непосредственно из определений.
\begin{lm}\label{homsmpr}
Пусть для $i=1,2$ заданы $\mathbin{\widehat{\otimes}}$-алгебра Хопфа $H_i$ и $H_i$-$\mathbin{\widehat{\otimes}}$-модульная алгебра~$A_i$.

\emph{(A)}~Тогда для гомоморфизма $\mathbin{\widehat{\otimes}}$-алгебр Хопфа $\psi\!:
H_1\to  H_2$ и  гомоморфизма $\mathbin{\widehat{\otimes}}$-алгебр $\varphi\!: A_1\to
A_2$, являющегося морфизмом $H_1$-$\mathbin{\widehat{\otimes}}$-модулей \emph{(}т.е.
$\varphi(h\cdot a)=\psi(h)\cdot\varphi(a)$\emph{)},  корректно определен
гомоморфизм  $\mathbin{\widehat{\otimes}}$-алгебр, заданный условием
\begin{equation}\label{morsmpr}
\varphi\mathop{\widehat{\#}} \psi\!:A_1\mathop{\widehat{\#}} H_1\to A_2\mathop{\widehat{\#}} H_2\!:a\otimes
h\mapsto \varphi(a)\otimes \psi(h).
\end{equation}

\emph{(B)}~Если, более того,  для $i=1,2$ $H_i$ кокоммутативна,  $A_i$  является $\mathbin{\widehat{\otimes}}$-алгеброй Хопфа и $H_i$-$\mathbin{\widehat{\otimes}}$-модульной биалгеброй, а $\varphi$  --- гомоморфизм  $\mathbin{\widehat{\otimes}}$-алгебр Хопфа, то
$\varphi\mathop{\widehat{\#}} \psi$ --- гомоморфизм  $\mathbin{\widehat{\otimes}}$-алгебр Хопфа.
\end{lm}

Рассмотрим диаграмму
\begin{equation}\label{diaanco}
   \xymatrix{
{\mathscr A}(G_1)\mathop{\widehat{\#}} {\mathscr A}(G_2)\ar[r]\ar[d]&{\mathscr A}(G_1\rtimes G_2) \ar[d]\\
{\mathscr A}_{\omega_1^\infty}(G_1)\mathop{\widehat{\#}} {\mathscr A}_{\omega_2^\infty}(G_2)\ar[r]&{\mathscr A}_{\omega^\infty}(G_1\rtimes G_2)\,, }
\end{equation}
где  верхняя горизонтальная стрелка --- изоморфизм из теоремы~\ref{AGsemid}, нижняяя горизонтальная стрелка --- изоморфизм из части~(C) теоремы~\ref{smprde}, левая вертикальная стрелка --- гомоморфизм  $\mathbin{\widehat{\otimes}}$-алгебр Хопфа из леммы~\ref{homsmpr}, а правая вертикальная стрелка --- гомоморфизм пополнения.
Несредственная проверка показывает, что выполнено следующее утверждение.
\begin{pr}\label{comdiaanco}
Диаграмма~\eqref{diaanco} коммутативна.
\end{pr}

\section{Теорема о разложении веса c экспоненциальным искривлением}
\label{s:TDWED}

Этот раздел содержит теорему об итерированном разложении веса c экспоненциальным искривлением на нильпотентной подгруппе. Хотя эта теорема и представляет самостоятельный интерес, для нас она является основным техническими результатом, который будут использован в \S\,\ref{s:DISPAM}.

Если  локально компактная группа $G$ компактно порождена (в частности, если она связна), то
\begin{equation}\label{wordwdef}
\omega(g)\!: = \min \{2^n \!: \, g \in U^{n} \},
\end{equation}
где $U$ --- относительно компактное подмножество, порождающее~$G$, и $U^0 =\{1\}$, является субмультипликативным весом. Более того,
его класс эквивалентности максимален относительно рассмотренного выше порядка (см. \cite[теорема~5.3]{Ak08} или \cite[Proposition~2.7]{ArAMN}). Поэтому всякий элемент этого класса мы будем называть \emph{максимальным субмультипликативным весом}\footnote{Буквальный перевод соответствующего английского термина \emph{word weight} для веса, определённого в~\eqref{wordwdef},  не благозвучен. Кроме того, для нас  свойство максимальности важнее явной конструкции.}. Изучение максимальных субмультипликативных весов является основным шагом при определении структуры $\widehat{\mathscr A}(G)$ для связных групп Ли (см. ниже).

\begin{df}
Пусть $\omega$ --- субмультипликативный вес на локально компактной группе $G$, а $\omega_H$ --- максимальный вес на её замкнутой подгруппе~$H$.  Будем говорить, что $\omega$ имеет \emph{экспоненциальное искривление
} на $H$, если
$$
\omega\preceq 1+\log \omega_H \quad\text{на $H$}.
$$
\end{df}
В силу того, что вес, определенный в~\eqref{wordwdef}, принимает значения в $[1,+\infty)$, 
здесь мы фактически переформулируем понятие экспоненциального искривления для функций длины, введённое в~\cite{ArUAPI}.
В частном случае, когда  $\omega$ --- максимальный вес на~$G$, мы получаем стандартное понятие подгруппы с экспоненциальным искривлением.

Так как интересующие нас вопросы о строении ${\widehat{{\mathscr A}}}(G)$ и
${\widehat{{\mathscr A}}}(G)^{\mathrm{PI}}$ для связной группы Ли легко сводятся к случаю, когда она линейна (см.~\S\,\ref{s:DISPAM}), в дальнейшем мы предполагаем, что $G$ --- связная линейная комплексная группа Ли.

Пусть $\mathfrak{g}$ --- её алгебра Ли. Напомним, что подгруппа, порождённая $\exp\mathfrak{h}$ для некоторой подалгебры Ли $\mathfrak{h}$ в $\mathfrak{g}$ называется \emph{интегральной}. Обозначим через $N$ \emph{нильпотентный  радикал}~$G$ --- пересечение ядер всех неприводимых конечномерных голоморфных представлений. Он является односвязной интегральной подгруппой, алгебра Ли которой есть $\mathfrak{n}$ --- нильпотентный  радикал алгебры Ли  $\mathfrak{g}$, т.е. пересечение ядер всех неприводимых конечномерных представлений, см., например, \cite[Proposition~1.1]{ArUAPI}.

Пусть $E$ обозначает \emph{экспоненциальный  радикал}~$G$. Мы опускаем его общее определение (термин  предложен в~\cite{Os02}, а в нужной нам общности оно дано в~\cite{Co08}; см. также \cite[\S\,3]{ArAnF}).
Для связных комплексных групп Ли мы используем как определение следующее свойство:
$E$ есть нормальная комплексная подгруппа Ли такая, что $G/E$ является наибольшей факторгруппой~$G$ локально изоморфной прямому произведению нильпотентной  и полупростой комплексных групп Ли  \cite[Proposition 3.11]{ArAnF}. (``Наибольшая'' она в том смысле, что каждая факторгруппа~$G$ с указанным свойством является факторгруппой~$G/E$.) Нам понадобятся только некоторые свойства экспоненциального  радикала.

\begin{pr}\label{exprpr} \cite[Proposition 1.3]{ArUAPI}
Для всякой связной линейной комплексной группы Ли её экспоненциальный  радикал $E$ является односвязной нильпотентной замкнутой нормальной подгруппой, содержащейся в нильпотентном  радикале $N$.
\end{pr}

Ниже мы предполагаем, что $N'$ ---  нормальная интегральная подгруппа в~$G$  такая, что $E\subset N'\subset N$. Отметим, что при указанных условиях $N'$ замкнута, нильпотентна и односвязна \cite[Lemma 1.2]{ArUAPI}. Также в \cite[Theorem 2.3]{ArUAPI} показано, что среди функций длины $G$, имеющих экспоненциальное искривление на $N'$ существует максимальная (с точностью до аддитивной(!) эквивалентности) и дано её асимптотическое разложение. Нам понадобится переформулировка этого результата в терминах весов.

Зафиксируем разложение $G=B\rtimes L$, где $B$ односвязна и разрешима, а $L$ линейно комплексно редуктивна (см. \cite[с.\,601, Theorem 16.3.7]{HiNe}). Заметим, что $N\subset B$  \cite[Proposition 4.44]{Le02}, а значит $N'\subset B$.
Пусть $\tau$ обозначает факторгомоморфизм $B\to B/N'$, $\mathfrak{b}$ и $\mathfrak{n}'$ --- алгебры Ли групп $B$ и $N'$ соответственно, $\mathfrak{h}$  --- некоторую подалгебру Картана в~$\mathfrak{b}$, а $\mathfrak{v}$ --- линейное подпространство дополняющее $\mathfrak{h}\cap \mathfrak{n}'$ в~$\mathfrak{h}$.
В силу \cite[Theorem 2.2]{ArUAPI}
$$
\mathfrak{n}'\times \mathfrak{v}\times L\to
G\!:(\eta,\xi,l)\mapsto \exp(\eta)\exp(\xi)\,l
$$
является биголоморфной эквивалентностью комплексных многообразий, а теорема~2.3 из \cite{ArUAPI} в терминах весов принимает следующий вид.

\begin{thm}\label{maindisto}
Пусть $G$ --- связная линейная комплексная группа Ли, $N'$ ---  её нормальная интегральная подгруппа  такая, что $E\subset N'\subset N$.

\emph{(A)}~Тогда существует вес $\omega_{\max}$ на $G$ с экспоненциальным искривлением на $N'$ такой, что $\omega\lesssim \omega_{\max}$ для всякого веса $\omega$ с экспоненциальным искривлением на $N'$.

\emph{(B)}~Пусть $\omega'$ и $\omega''$ --- некоторые максимальные субмультипликативные веса на~$B/N'$ и~$L$ соответственно. Тогда $\omega_{\max}$ эквивалентен функции
\begin{equation}\label{thredecw}
\upsilon(g)\!:=(1+\|\eta\|)\;\omega'\tau(\exp(\xi))\;\omega''(l),
\end{equation}
где $g$ единственным образом представлен в виде $\exp(\eta)\exp(\xi)\,l$, $\eta\in\mathfrak{n}'$, $\xi\in\mathfrak{v}$ и $l\in L$, а $\|\cdot\|$ --- некоторая норма на $\mathfrak{n}'$ \emph{(}как на векторном пространстве\emph{)}.
\end{thm}

Назовём $\omega_{\max}$ \emph{максимальным весом с экспоненциальным искривлением на $N'$}, ср. \cite[Definition 2.4]{ArUAPI}.

Нас, в первую очередь, интересуют крайние случаи: $N'=E$ и $N'=N$, когда ${\mathscr A}_{\omega^\infty}$ совпадает с ${\widehat{{\mathscr A}}}(G)$ и ${\widehat{{\mathscr A}}}(G)^{\mathrm{PI}}$ (см. теорему~\ref{redtog}). Однако естественно доказать теорему о разложении для произвольной~$N'$.

Здесь мы покажем, что  из структурной теории линейных групп Ли и теоремы~\ref{maindisto}  следует, что максимальный вес с экспоненциальным искривлением на~$N'$ может быть представлен  в виде, согласованном с разложением~$G$ в некоторое итерированное полупрямое произведение, т.е. на каждом шаге выполнено условие~\eqref{eqom1om2}, гарантирующее разложение в смэш-произведение. Более того, веса на соответствующих множителях имеют довольно простой вид.

\begin{thm}\label{itsemdd0}
Пусть $G$ --- связная линейная комплексная группа Ли, $N'$ ---  её нормальная интегральная подгруппа  такая, что $E\subset N'\subset N$, а $\omega_{\max}$ --- субмультипликативный вес на~$G$, максимальный среди весов с экспоненциальным искривлением на~$N'$. Зафиксируем разложение $G=B\rtimes L$, где $B$ односвязна и разрешима, а $L$ линейно комплексно редуктивна. Тогда $G$ может быть разложена в итерированное полупрямое произведение
\begin{equation}\label{itdirp0}
 G=((\cdots (F_1 \rtimes F_2)\rtimes\cdots)
\rtimes F_n),
\end{equation}
где $F_1\cong\cdots\cong F_{n-1}\cong\mathbb{C}$ и
$F_n=L$, с выполнением следующих условий. Положим $G_1\!:=F_1$ и $G_i\!:=G_{i-1}\rtimes F_i$ при $i\le 2$, тогда $G_p=N'$ (где $p$ --- размерность $N'$) и $G_{n-1}=B$. Кроме того, для каждого $i$ существует субмультипликативный вес~$\omega_i$ на~$F_i$, такой что выполнено следующее.

\emph{(A)}
\begin{equation}\label{omiz}
\omega_i(z)=\left\{
\begin{array}{ll}
1+|z|,&\quad \text{если $i\le p$,}\\
 \exp(|z|^{1/w_{n-i}}),&\quad \text{если $p<i\le n-1$,}
 \end{array}
 \right.\qquad(z\in \mathbb{C}),
\end{equation}
где $w_1,\ldots, w_{n-1-p}$ --- неубывающая последовательность в $\mathbb{N}$ c
$w_1=1$, а $\omega_n$ --- максимальный субмультипликативный вес на~$L$.

\emph{(B)} Имеет место итерированное разложение $\omega_{\max}$ в следующем смысле. Пусть $\widetilde\omega_i$ --- ограничение~$\omega_{\max}$ на~$G_i$. Тогда $\widetilde\omega_1 \simeq  \omega_1$ и для всех $i=1,\ldots,n-1$
$$
\widetilde\omega_{i+1}(gf)\simeq  \widetilde\omega_i(g)\omega_{i+1}(f)\quad\text{на
$G_i\times F_{i+1}$} \quad(\text{$g\in G_i$, $f\in F_{i+1}$}),
$$
в частности,
$\omega_{\max}(gf)\simeq  \widetilde\omega_{n-1}(g)\omega_n(f)$.
\end{thm}
\begin{proof}
В силу теоремы~\ref{maindisto}
\begin{equation}\label{wlfsec0}
\omega_{\max}(\exp(\eta)\exp(\xi)l)\simeq (1+\|\eta\|)\;\omega'(\tau(\exp(\xi)))\;\omega''(l).
\end{equation}
где $\|\cdot\|$ --- некоторая норма на $\mathfrak{n}'$ (как на векторном пространстве).
Используя это асимпотическое разложение, мы сначала построим разложение~$N'$ в итерированное полупрямое произведение
c соответствующими субмультипликативными весами, потом продолжим его на~$B$, а затем на $G$.

Можно предполагать, что подгруппа~$N'$ нетривиальна. Будучи подалгеброй нильпотентной алгебры Ли $\mathfrak{n}$,
алгебра $\mathfrak{n}'$ также нильпотентна. Поэтому она может быть представлена как итерированная полупрямая сумма
одномерных алгебр Ли:
$$
\mathfrak{n}'\cong ((\cdots (\mathfrak{f}_1 \rtimes \mathfrak{f}_2)\rtimes\cdots)
\rtimes\mathfrak{f}_p),
$$
где $\mathfrak{f}_i\!:=\mathbb{C} e_i$ для некоторого $e_i\in \mathfrak{n}'$, а
$p$ --- размерность $N'$. Рассмотрим
подгруппы $F_i\!:=\exp \mathfrak{f}_i$, $i=1,\ldots, p$, в $G$.
Тогда $N'$ представляется как итерированное полупрямое произведение:
$$
N'\cong ((\cdots (F_1 \rtimes F_2)\rtimes\cdots) \rtimes F_p).
$$

Определим $G_1,\ldots,G_p$ как соответствующие итерированные полупрямые произведения и рассмотрим для каждого
$i$ ограничение $\widetilde\omega_i$  веса $\omega_{max}$ на~$G_i$.
Возьмём в качестве $\|\cdot\|$ норму $\sum_i|s_i|$, где
$(s_1,\ldots,s_p)$ --- координаты в базисе $e_1,\ldots,e_p$. Так как
$$
1+\sum_i|s_i|\le \prod_i(1+|s_i|)\le \Bigl(1+\sum_i|s_i|\Bigr)^p,
$$
получаем
\begin{equation}\label{logprodeq0}
1+\|\eta\|\simeq \prod_i (1+|s_i|)\qquad \bigl(\eta=\sum_i s_i e_i\bigr)\,.
\end{equation}
Из \eqref{wlfsec0} и \eqref{logprodeq0} следует, что для $i=1,\ldots, p$
$$
\widetilde\omega_{i+1}(gz)\simeq \widetilde\omega_i(g)\,(1+|z|)\quad \text{на $G_i\times
F_{i+1}$},
$$
и мы можем положить $\omega_i(z)\!:=1+|z|$.

Теперь продолжим разложение на~$B$. Заметим, что $B/N'$ также односвязна и, более того, нильпотентна, так как таковой является группа $B/E$  \cite[\S\,3]{ArAnF} (напомним, что $E\subset N'$ согласно предположению).
Обозначим через $\mathfrak{q}$ алгебру Ли группы $B/N'$, а через ${\mathscr F}$ --- её нижний центральный ряд (он определяется индуктивно как $\mathfrak{q}_1\!:=\mathfrak{q}$, $\mathfrak{q}_k\!:=[\mathfrak{q}, \mathfrak{q}_{k-1}]$). Зафиксируем ${\mathscr F}$-базис $y_1,\ldots,y_{n-1-p}$ в~$\mathfrak{q}$; это означает, что
последовательность
\begin{equation}\label{wkdef}
w_k\!:=\max\{j:\,y_k\in \mathfrak{q}_j \}
\end{equation}
не убывает и $\mathfrak{q}_j$ есть линейная оболочка $\{y_k : w_k \ge j\}$ для каждого~$j$. (Мы обозначаем размерность~$B$ и~$\mathfrak{b}$ через $n-1$).

Пусть $\mathfrak{h}$ и $\mathfrak{v}$ выбраны как выше, т.е.~$\mathfrak{h}$ --- некоторая подалгебра Картана алгебры~$\mathfrak{b}$, а подпространство $\mathfrak{v}$  --- некоторое линейное дополнение подпространства $\mathfrak{h}\cap \mathfrak{n}'$ в $\mathfrak{h}$.   Обозначим через $\mathfrak{e}$  алгебру Ли группы $E$. Тогда $\mathfrak{b}=\mathfrak{e}+\mathfrak{h}$, см. \cite[Theorem 2.2]{ArUAPI}. Так как $\mathfrak{e}\subset \mathfrak{n}'$, имеем $\mathfrak{b}\subset \mathfrak{n}'+\mathfrak{h}\subset \mathfrak{b}$. Отсюда и из определения $\mathfrak{v}$ вытекает, что $\mathfrak{b}=\mathfrak{n}'\oplus \mathfrak{v}$.

Итак, факторотображение алгебр Ли $\mathfrak{b}\to \mathfrak{q}$ отображает $\mathfrak{v}$ биективно на~$\mathfrak{q}$.
Следовательно, прообразы $e_{p+1},\ldots,e_{n-1}$  элементов $y_{n-1-p},\ldots, y_1$ составляют базис в~$\mathfrak{v}$ (мы обращаем порядок индексов). Так как $y_1,\ldots,y_{n-1-p}$ является ${\mathscr F}$-базисом, мы можем записать $\mathfrak{b}$ как итерированную полупрямую сумму:
$$
\mathfrak{b}\cong ((\cdots (\mathfrak{n}' \rtimes \mathfrak{f}_{p+1})\rtimes\cdots)
\rtimes\mathfrak{f}_{n-1})\,,
$$
где $\mathfrak{f}_{p+1}\!:=\mathbb{C} e_{p+1},\ldots,\mathfrak{f}_{n-1}\!:=\mathbb{C} e_{n-1}$.
Рассмотрим подгруппы $F_i\!:=\exp \mathfrak{f}_i$  в $G$ для $i=p+1,\ldots,
n-1$. Тогда
$$
B\cong ((\cdots (N' \rtimes F_{p+1})\rtimes\cdots) \rtimes
F_{n-1}).
$$
Определим $G_{p+1},\ldots,G_{n-1}$ как соответствующие итерированные полупрямые произведения и рассмотрим для каждого
$i$ ограничение $\widetilde\omega_i$  веса $\omega_{max}$ на~$G_i$.

Так как $B/N'$ односвязна и нильпотентна, мы можем применить теорему~3.1 из
\cite{ArAMN} (в формулировке для весов), согласно которой $\omega'$ (максимальный вес на $B/N'$) эквивалентен  функции
\begin{equation*}
 \mu(q):=\exp(\max_k|s_k|^{1/w_k}),\qquad (q\in B/N')
\end{equation*}
где $(s_1,\ldots, s_{n-1-p})$ --- канонические координаты второго рода на $B/N'$, ассоциированные с базисом $y_1,\ldots,y_{n-1-p}$.

Очевидно, $\mu$ эквивалента функции  $q\mapsto \prod_k\exp(|s_k|^{1/w_k})$. Тогда из \eqref{wlfsec0} следует, что
$$
\widetilde\omega_{i+1}(gz)\simeq \widetilde\omega_i(g)\,\exp(|z|^{1/w_{n-i}})\quad \text{на $G_i\times F_{i+1}$}\qquad(i\ge p).
$$
Тем самым мы можем положить $\omega_i(z)\!:=\exp(|z|^{1/w_{n-i}})$.

Итак, мы построили разложение для $B=G_{n-1}$. Осталось сделать последний шаг.
Положим $F_n\!:=L$ и $\omega_n\!:=\omega''$ (максимальный вес на~$L$). Тогда, снова применяя~\eqref{wlfsec0}, получаем, что
$$
\widetilde\omega_n(gf)\simeq \widetilde\omega_{n-1}(g)\,\omega_n(f)\quad \text{на $G_{n-1}\times F_n$.}
$$
Так как $\widetilde\omega_n= \omega_{\max}$, это завершает доказательство.
\end{proof}
\begin{rem}\label{seqwi}
В приведённом выше доказательстве теоремы~\ref{itsemdd0} мы получили немного больше информации, чем указано в формулировке. А именно, в~\eqref{wkdef} выписан явный вид последовательности $w_1,\ldots, w_{n-1-p}$. Обозначим через~$m$ класс нильпотентности $\mathfrak{q}$, т.е. $\mathfrak{q}_{m+1}=0$, но $\mathfrak{q}_{m}\ne 0$. Тогда из указанной формулы следует, что $w_{n-1-p}=m$. Более того, если $i\in\{1,\ldots,m\}$, то количество членов последовательнсоти $w_1,\ldots, w_{n-1-p}$, равных~$j$, совпадает с размерностью $\mathfrak{q}_j/\mathfrak{q}_{j+1}$.
\end{rem}

\section{Универсальные алгебры}
\label{s:UNAL}

В этом разделе мы обсуждаем универсальные свойства интересующих нас пополнений алгебры ${\mathscr A}(G)$. Эти свойства будут использованы в \S\,\ref{s:DISPAM} для исследования крайних случаев, $N'=E$ и $N'=N$.  Результаты из \S\,\ref{s:TDWED}  здесь не понадобятся.

\subsection*{Общее универсальное свойство}

Сначала мы покажем, что для всякого максимального веса $\omega_{max}$  с экспоненциальным искривлением на~$N'$ алгебра ${\mathscr A}_{\omega_{max}^\infty}(G)$ обладает важным универсальным свойством. Предварительно напомним, что всякий голоморфный гомоморфизм $\pi$ из комплексной группы Ли $G$ в группу обратимых элементов ${\mathop{\mathrm{GL}}\nolimits}(A)$ банаховой алгебры $A$ индуцирует непрерывный гомоморфизм $\mathbin{\widehat{\otimes}}$-алгебр $\bar\pi\!:{\mathscr A}(G)\to A$.
При этом выполнено соотношение
\begin{equation}\label{pixg3}
\pi(g)= \bar\pi(\delta_g) \qquad (g\in G).
\end{equation}
Более того, соответствие $\pi\mapsto \bar\pi$ биективно (\cite{Li70}, см. также \cite[\S\,5]{ArAnF}).

Заметим также, что для голоморфного гомоморфизма $\pi\!:G\to {\mathop{\mathrm{GL}}\nolimits}(A)$, где  $A$ --- банахова алгебра с нормой $\|\cdot\|$, функция $g\mapsto\|\pi(g)\|$ является субмультипликативным весом.

\begin{thm}\label{univdes0}
Пусть $G$ --- связная линейная комплексная группа Ли, а $N'$ ---  её нормальная интегральная подгруппа  такая, что $E\subset N'\subset N$. Пусть также $\omega_{max}$ --- субмультипликативный вес на~$G$, максимальный среди весов с экспоненциальным искривлением на $N'$, а $\iota\!:{\mathscr A}(G)\to {\mathscr A}_{\omega_{max}^\infty}(G)$ --- гомоморфизм пополнения.
Тогда для всякого голоморфного гомоморфизма $\pi\!:G\to {\mathop{\mathrm{GL}}\nolimits}(A)$, где  $A$ --- банахова алгебра с нормой $\|\cdot\|$, такого что вес $g\mapsto\|\pi(g)\|$ имеет экспоненциальное искривление на~$N'$, существует единственный непрерывный гомоморфизм
$\rho\!:{\mathscr A}_{\omega_{max}^\infty}(G) \to A$, делающий диаграмму
\begin{equation*}
  \xymatrix{
{\mathscr A}(G)  \ar[r]^{\iota}\ar[rd]_{\bar\pi}& {\mathscr A}_{\omega_{max}^\infty}(G)  \ar@{-->}[d]^{\rho}\\
 &A\\
 }
\end{equation*}
коммутативной.
\end{thm}

Для доказательства понадобится следующее утверждение.

\begin{pr}\label{nornom}
Пусть $G$ --- комплексная группа Ли.

\emph{(A)}~Если $\|\cdot\|$ --- непрерывная субмультипликативная преднорма на ${\mathscr A}(G)$ и $\omega(g)\!:=\|\delta_g\|$, то $\|x\|\le \|x\|_\omega$ для каждого $x\in{\mathscr A}(G)$ \emph{(}определение $\|\cdot\|_\omega$  см. в~\eqref{Vup}\emph{)}.

\emph{(B)}~Если $A$ --- банахова алгебра с нормой $\|\cdot\|$, $\pi\!:G\to {\mathop{\mathrm{GL}}\nolimits}(A)$ --- голоморфный гомоморфизм, а $\omega(g)\!:=\|\pi(g)\|$, то $\bar\pi\!:{\mathscr A}(G)\to A$ продолжается до непрерывного гомоморфизма ${\mathscr A}_{\omega}(G)\to A$.
\end{pr}
\begin{proof}
(A)~Утверждение следует из \cite[теорема 5.2(1)]{Ak08} (обозначения в её формулировке см. там же, \S\,3.4.1).

(B)~Так как согласно~\eqref{pixg3} $\pi(g)= \bar\pi(\delta_g)$ для всех $g\in G$, мы можем применить утверждение~(A) к $\omega$. Тем самым $\|\bar\pi(x)\|\le \|x\|_{\omega}$ для каждого $x\in{\mathscr A}(G)$ и, таким образом, мы можем продолжить $\bar\pi$ до непрерывного гомоморфизма ${\mathscr A}_{\omega}(G)\to A$.
\end{proof}

\begin{proof}[Доказательство теоремы~\ref{univdes0}]
Предположим, что вес $\omega'(g)\!:=\|\pi(g)\|$ имеет экспоненциальное искривление на~$N'$.
Так как $\omega_{max}$ ---  максимальный среди таких весов, то существуют $C>0$ и $n\in\mathbb{N}$, такие что $\omega'(g)\le C\,\omega_{max}(g)^n$ для всех $g\in G$. Отсюда следует, что тождественное отображение на ${\mathscr A}(G)$ продолжается до непрерывного гомоморфизма банаховых алгебр ${\mathscr A}_{\omega_{max}^n}(G)\to{\mathscr A}(G)_{\omega'}$.
С  другой стороны, в силу части~(B) предложения~\ref{nornom} мы можем продолжить $\bar\pi$ до непрерывного гомоморфизма ${\mathscr A}_{\omega'}(G)\to A$. Дальше ясно.
\end{proof}

\subsection*{Оболочки относительно класса банаховых PI-алгебр}

Напомним, что ассоциативная алгебра $A$ (в нашем случае над $\mathbb{C}$) называется PI-\emph{алгеброй}, если она удовлетворяет полиномиальному тождеству, т.е. существуют $n\in\mathbb{N}$ и ненулевой элемент $p$ свободной алгебры с $n$ образующими, такие что  $p(a_1,\ldots,a_n)=0$ для всех $a_1,\ldots,a_n\in A$. Общие сведения о PI-алгебрах см. например, \cite{KKR16}.

Если группа линейна, то она, по определению, допускает инъективное голоморфное конечномерное представление, которое можно трактовать как гомоморфизм в конечномерную алгебру. Так как всякая конечномерная алгебра является PI-алгеброй, класс банаховых PI-алгебр естественно появляется в наших рассмотрениях.

\begin{rem}
В этой статье нет ссылок непосредственно на результаты из алгебраической теории полиномиальных тождеств.
Единственное утверждение о PI-алгебрах, которое нам понадобится, --- это теорема~\ref{AGcrit}, доказанная ниже. Она связывает веса с экспоненциальным искривлением и голоморфные гомоморфизмы из $G$ в банаховы PI-алгебры. Это утверждение является уточнением теоремы~1.14 из~\cite{ArUAPI}, которая, в свою очередь, получается переформулировкой для групп Ли теоремы~2 из~\cite{ArPiLie}, доказанной для алгебр Ли. А вот последний результат уже существенно опирается на (довольно глубокие) алгебраические результаты о PI-алгебрах.

Следует добавить, что теорема~\ref{maindisto}, приведённая выше, (точнее, её вариант для функций длины --- \cite[Theorem 2.3]{ArUAPI}) хотя и не содержит в формулировке упоминаний PI-алгебр, также
получен с использованием теоремы~1.14 из~\cite{ArUAPI}.
\end{rem}

Следуя \cite[Chapter~I, \S\,1, с.\,4, Definition 1.1]{Ma86}, мы называем ассоциативную алгебру, снабжённую  структурой топологического векторного пространства, такой что умножение раздельно непрерывно, \emph{топологической алгеброй}. Ниже будут использованы проективные пределы в категории топологических алгебр, см. подробности в \cite[Chapter~III]{Ma86}.

\begin{df}
Будем говорить, что алгебра Аренса-Майкла \emph{локально содержится в классе банаховых} PI-\emph{алгебр}, если она является проективным пределом банаховых PI-алгебр в категории топологических алгебр.
\end{df}
Отметим, что пределы проективной системы банаховых алгебр в категориях топологических алгебр и алгебр Аренса-Майкла совпадают.

Проективные пределы банаховых PI-алгебр встречаются в некоммутативной ком\-плекс\-но-ана\-ли\-ти\-ческой геометрии (см., например, результат Люминета от том, что рассмотренные  Дж.Тэйлором алгебра свободных голоморфных функций и её обобщения принадлежат этому классу \cite[Proposition 2.6]{Lu86}).

\begin{df}\label{BPIun}
Пусть $A$ --- ассоциативная топологическая алгебра. Мы говорим, что пара, состоящая из  алгебры Аренса-Майкла $\widehat A^{\,\mathrm{PI}}$, которая локально содержится в классе банаховых PI-алгебр, и непрерывного гомоморфизма $\iota\!:A\to\widehat A^{\,\mathrm{PI}}$, является \emph{оболочкой $A$ относительно  класса банаховых} PI-\emph{алгебр}, если выполнено следующее универсальное свойство. Для всякой банаховой PI-алгебры~$B$ и всякого непрерывного гомоморфизма
$\varphi\!: A\to B$ существует единственный непрерывный гомоморфизм
$\widehat\varphi\!:\widehat A^{\,\mathrm{PI}} \to B$, такой что диаграмма
\begin{equation}\label{AMen}
  \xymatrix{
A \ar[r]^{\iota}\ar[rd]_{\varphi}&\widehat A^{\,\mathrm{PI}}\ar@{-->}[d]^{\widehat\varphi}\\
 &B\\
 }
\end{equation}
коммутативна.
\end{df}

\begin{lm}\label{PIeneqd}
Пусть $A$ --- ассоциативная топологическая алгебра. Пара $(\widehat A^{\,\mathrm{PI}},\iota)$ является оболочкой~$A$ относительно  класса банаховых \emph{PI}-алгебр тогда и только тогда, когда для всякой~$B$, локально содержащейся в классе банаховых \emph{PI}-алгебр, и всякого непрерывного гомоморфизма
$\varphi\!: A\to B$ существует единственный непрерывный гомоморфизм
$\widehat\varphi\!:\widehat A^{\,\mathrm{PI}} \to B$, такой что диаграмма~\eqref{AMen} коммутативна.
\end{lm}
\begin{proof}
Достаточность очевидна. Необходимость следует из того, что согласно определению $B$ является
проективным пределом банаховых PI-алгебр и универсального свойства проективного предела.
\end{proof}

\begin{pr}
Для всякой ассоциативной топологической алгебры оболочка относительно  класса банаховых \emph{PI}-алгебр существует и единственна с точностью до изоморфизма.
\end{pr}
\begin{proof}
Рассмотрим ассоциативную топологическую алгебру~$A$.
Пусть $\widehat A^{\,\mathrm{PI}}$ ---  пополнение $A$ относительно топологии,  заданной семейством всех непрерывных субмультипликативных преднорм, для которых соответствующие пополнения есть (банаховы) PI-алгебры, а $\iota$ --- гомоморфизм пополнения.  Это доказывает существование. Единственность следует из единственности $\widehat\varphi$ в лемме~\ref{PIeneqd}.
\end{proof}

Легко видеть из определения, что соответствие $A\mapsto \widehat A^{\,\mathrm{PI}}$ продолжается до функтора из категории ассоциативных топологических алгебр в полную подкатегорию категории алгебр Аренса-Майкла, объектами которой являются  алгебры, локально содержащиеся в классе банаховых PI-алгебр.

\subsection*{Оболочки и экспоненциальное искривление}
Далее нас интересуют две оболочки алгебры ${\mathscr A}(G)$, а именно, ${\widehat{{\mathscr A}}}(G)^{\mathrm{PI}}$ (оболочка относительно  класса банаховых PI-алгебр) и ${\widehat{{\mathscr A}}}(G)$ (оболочка относительно  класса всех банаховых алгебр, иначе говоря, оболочка Аренса-Майкла, см. конец \S\,\ref{s:ASP}).

Для произвольной комплексной группы Ли $G$ обозначим через $\Lambda$ её линеаризатор, т.е. пересечение ядер всевозможных конечномерных голоморфных представлений. Так как $\Lambda$ --- нормальная подгруппа, то можно рассмотреть отображение $G\to G/\Lambda$, которое индуцирует гомоморфизмы
$$
\alpha\!:{\widehat{{\mathscr A}}}(G)\to {\widehat{{\mathscr A}}}(G/\Lambda)\quad\text{и}\quad
\beta\!:{\widehat{{\mathscr A}}}(G)^{\mathrm{PI}}\to {\widehat{{\mathscr A}}}(G/\Lambda)^{\mathrm{PI}}.
$$

\begin{pr}\label{isoLa}
Если $G$ связна, то $\alpha$ и $\beta$ являются топологическими изоморфизмами.
\end{pr}
\begin{proof}
Утверждение для $\alpha$ доказано в части~(A) теоремы~5.3 из \cite{ArAnF}. (Рассуждение использует
теорему~2.2 из \cite{ArAnF}, которая утверждает что для связных групп~$\Lambda$ совпадает с обобщённым линеаризатором, т.е. пересечением ядер всевозможных голоморфных гомоморфизмов в группы обратимых элементов банаховых алгебр.) Утверждение для $\beta$ доказывается точно так же с тем отличием, что вместо класса всех банаховых алгебр надо использовать класс банаховых PI-алгебр.
\end{proof}

Как хорошо известно,   $G/\Lambda$ линейна (т.е. она является линеаризацией $G$), поэтому предложение~\ref{isoLa} сводит вопросы о строении ${\widehat{{\mathscr A}}}(G)$ и ${\widehat{{\mathscr A}}}(G)^{\mathrm{PI}}$ к случаю, когда $G$ линейна.

Следующая теорема даёт критерий того, что пополнение ${\mathscr A}(G)$ является банаховой PI-алгеброй, в терминах экспоненциального искривления.

\begin{thm}\label{AGcrit}
Пусть $G$ --- связная линейная комплексная группа Ли, а $N$ --- её нильпотентный радикал. Предположим, что $A$ --- банахова алгебра с нормой $\|\cdot\|$, а $\pi\!:G\to {\mathop{\mathrm{GL}}\nolimits}(A)$ --- голоморфный гомоморфизм, такой что ассоциированный гомоморфизм $\bar\pi\!:{\mathscr A}(G)\to A$ \emph{(}см. \eqref{pixg3}\emph{)} имеет плотный образ. Тогда $A$ является \emph{PI}-алгеброй, если и только если вес  $g\mapsto\|\pi(g)\|$ имеет экспоненциальное искривление на $N$.
\end{thm}

Чтобы доказать теорему, мы используем два отображения и лемму об их свойствах.

Для произвольной комплексной группы Ли $G$ рассмотрим гомоморфизм
\begin{equation}\label{taudef}
\sigma\!:U(\mathfrak{g}) \to \mathscr{A}(G)\!:\langle \sigma(X), f\rangle\!:= [{\widetilde X}f](1) \qquad (X\in U(\mathfrak{g}),\, f \in  \mathcal{O}(G)),
\end{equation}
где ${\widetilde X}$ --- левоинвариантный дифференциальный оператор, соответствующий $X\in U(\mathfrak{g})$.

Группу обратимых элементов ${\mathop{\mathrm{GL}}\nolimits}(A)$ банаховой алгебры $A$ можно рассмотреть как комплексную группу Ли (возможно бесконечномерную), а голоморфный гомоморфизм $\pi\!:G\to {\mathop{\mathrm{GL}}\nolimits}(A)$ из  комплексной группы Ли $G$ --- как гомоморфизм групп Ли. При этом алгебра Ли группы ${\mathop{\mathrm{GL}}\nolimits}(A)$ отождествляется с $A$, а экспоненциальное отображение с
$$
\exp\!:A\to {\mathop{\mathrm{GL}}\nolimits}(A)\!:a\mapsto \sum_{n=0}^\infty \frac{a^n}{n!},
$$
см., например, \cite[Example III.1.11(b) и Remark IV.2.2]{Ne05}.
Обозначим, как и выше, алгебру Ли группы $G$ через $\mathfrak{g}$.
Применяя к $\pi$ функтор Ли, получаем гомоморфизм банаховых алгебр Ли $\mathrm{L}_\pi\!:\mathfrak{g}\to A$, такой что $\pi\exp=\exp\mathrm{L}_\pi$.

\begin{lm}\label{dersi}
Пусть $G$ --- комплексная группа Ли, $A$ --- банахова алгебра, $\pi\!:G\to {\mathop{\mathrm{GL}}\nolimits}(A)$ --- голоморфный гомоморфизм, а
$\theta\!:U(\mathfrak{g})\to A$ --- гомоморфизм, порождённый гомоморфизмом алгебр Ли $\mathrm{L}_\pi\!: \mathfrak{g} \to A$. Тогда
$\theta=\bar\pi\sigma$.
\end{lm}
\begin{proof}
Достаточно показать, что $\bar\pi\sigma(X)=\mathrm{L}_\pi(X)$ для каждого $X\in\mathfrak{g}$.

Зафиксируем $X$. В силу определения экспоненциального отображения имеем
$$[{\widetilde X}f(g)]=\frac{d}{dz}\Bigl |_{z=0}f(g\exp(zX))\qquad(z\in\mathbb{C})$$
для всех $f\in \mathcal{O}(G)$ и $g\in G$.
Напомним, что $\bar\pi$ задаётся формулой
\begin{equation}\label{pixg2}
\langle F, \bar\pi(a')\rangle =\langle a', \pi_F\rangle,
\end{equation}
где $a'\in {{\mathscr A}}(G)$, $F$ --- непрерывный линейный функционал на $A$, а
функция $\pi_F\in \mathcal{O}(G)$ определена соотношением
\begin{equation}\label{pixg4}
\pi_F(g)\!:= \langle F, \pi(g)\rangle\qquad(g\in G),
\end{equation}
см. \cite{Li70,Li72}.
Тем самым, учитывая~\eqref{taudef}, для каждого $F$ имеем
$$
\langle F, \bar\pi\sigma(X)\rangle =\langle \sigma(X), \pi_F\rangle=[{\widetilde X}\pi_F](1)=\frac{d}{dz}\Bigl |_{z=0}\pi_F(\exp(zX))=\frac{d}{dz}\Bigl |_{z=0} \langle F, \pi(\exp(zX))\rangle.
$$
Так как $\pi\exp=\exp\mathrm{L}_\pi$, то
$$
\frac{d}{dz}\Bigl |_{z=0}\pi(\exp(zX))=\frac{d}{dz}\Bigl |_{z=0}\exp(\mathrm{L}_\pi(zX))=\mathrm{L}_\pi(X).
$$

Итак, $\langle F, \bar\pi\sigma(X)\rangle =\langle F, \mathrm{L}_\pi(X)\rangle$.  В силу произвольности $F$ заключаем, что $\bar\pi\sigma(X)=\mathrm{L}_\pi(X)$, что и требовалось.
\end{proof}

\begin{proof}[Доказательство теоремы~\ref{AGcrit}]
Если гомоморфизм $\theta\!:U(\mathfrak{g})\to A$, порождённый $\mathrm{L}_\pi\!: \mathfrak{g} \to A$, имеет плотный образ, то $A$ является PI-алгеброй в том и только том, случае, когда $g\mapsto\|\pi(g)\|$ имеет экспоненциальное искривление на $N$ \cite[Theorem 1.14]{ArUAPI}\footnote{Напомним, что отношение $\preceq$ в  \cite{ArUAPI} получается из используемого здесь логарифмированием.}.

Заметим, что $\sigma\!:U(\mathfrak{g})\to {\mathscr A}(G)$ имеет плотный образ. Это следует из того, что соответствующее линейное отображение $\sigma'\!:\mathcal{O}(G)\to U(\mathfrak{g})'$  сильных дуальных пространств инъективно для каждой связной группы Ли, см., например, обсуждение  после формулы~(42) в~\cite{Pir_stbflat}. Кроме того, образ $\bar\pi$ плотен согласно предположению, а  $\theta=\bar\pi\sigma$ в силу леммы~\ref{dersi}. Следовательно, образ $\theta$ также плотен, и мы можем применить упомянутый выше результат из \cite{ArUAPI}.
\end{proof}

\begin{thm}\label{exBaPI}
Пусть $G$ --- связная линейная комплексная группа Ли.
Если $\omega$ --- субмультипликативный вес на~$G$ c экспоненциальным искривлением на $N$, то ${\mathscr A}_{\omega}(G)$ является \emph{PI}-алгеброй.
\end{thm}
\begin{proof}
По определению, банахова алгебра ${\mathscr A}_{\omega}(G)$ является пополнением ${\mathscr A}(G)$ относительно
преднормы $\|\cdot\|_\omega$, о последней см.~\eqref{Vup} и ниже. Обозначим через $\bar\pi$ гомоморфизм пополнения ${\mathscr A}(G)\to{\mathscr A}_{\omega}(G)$, а через $\pi$ --- соответствующий ему голоморфный гомоморфизм  $G\to{\mathop{\mathrm{GL}}\nolimits}({\mathscr A}_{\omega}(G))$, см.~\eqref{pixg3} и дальнейшие ссылки. Так как $\pi(g)= \bar\pi(\delta_g)$ согласно упомянутой формуле~\eqref{pixg3}, то из~\eqref{Vup} следует, что $\|\pi(g)\|_\omega\le \omega(g)$. Поскольку $\omega$ имеет экспоненциальное искривление на~$N$, отсюда следует, что $g\mapsto \|\pi(g)\|_\omega$ также имеет экспоненциальное искривление на~$N$. Очевидно, образ $\bar\pi$ плотен. Поэтому ${\mathscr A}_{\omega}(G)$ является PI-алгеброй согласно теореме~\ref{AGcrit}.
\end{proof}

Теперь мы можем доказать основной результат этого параграфа.

\begin{thm}\label{redtog}
Пусть $G$ --- связная линейная комплексная группа Ли.

\emph{(A)}~Если $\omega$ --- субмультипликативный вес на~$G$, максимальный среди весов с экспоненциальным искривлением на~$E$, то он  является максимальным субмультипликативным весом и ${\mathscr A}_{\omega^\infty}(G)\cong \widehat{\mathscr A}(G)$.

\emph{(B)}~Если $\omega$ --- субмультипликативный вес на~$G$, максимальный среди весов с экспоненциальным искривлением на $N$, то ${\mathscr A}_{\omega^\infty}(G)\cong \widehat{\mathscr A}(G)^{\mathrm{PI}}$.
\end{thm}

\begin{proof}
(A)~Первое утверждение фактически доказано в \cite{Co08}. Действительно, экспоненциальный радикал связной группы Ли является подгруппой с экспоненциальным искривлением \cite[Theorem~6.3]{Co08}. В частности, всякий субмультипликативный  вес  имеет экспоненциальное искривление на $E$. Следовательно, $\omega$ является максимальным в классе всех субмультипликативных весов.

Отсюда получаем, что ${\mathscr A}_{\omega^\infty}(G)\cong \widehat{\mathscr A}(G)$ (в силу \cite[теоремы~5.3 и~6.2]{Ak08}, см. также   \cite[Proposition 2.8]{ArAMN}).

(B)~Если $n\in\mathbb{N}$, то $\omega^n\simeq \omega$ и из теоремы~\ref{exBaPI} следует, что ${\mathscr A}_{\omega^n}(G)$ является PI-алгеброй.
Таким образом, ${\mathscr A}_{\omega^\infty}(G)$  локально содержится в классе банаховых PI-алгебр.

Осталось показать, что ${\mathscr A}_{\omega^\infty}(G)$ обладает универсальным свойством из определения~\ref{BPIun}.
Пусть $A$ --- банахова PI-алгебра и $\bar\pi\!:{\mathscr A}(G)\to A$ --- непрерывный гомоморфизм. Можно предполагать, что $\bar\pi$ имеет плотный образ. В силу теоремы~\ref{AGcrit} вес  $\omega_1(g)\!:=\|\pi(g)\|$ имеет экспоненциальное искривление на $N$.
Из части~(B) предложения~\ref{nornom} следует, что $\bar\pi$  продолжается до непрерывного гомоморфизма ${\mathscr A}_{\omega_1}(G)\to A$. Так как $\omega$ --- субмультипликативный вес на~$G$, максимальный среди с весов экспоненциальным искривлением на $N$, то $\omega_1 \preceq\omega$. Тем самым найдётся $n\in\mathbb{N}$ такой, что $\bar\pi$  продолжается до непрерывного гомоморфизма ${\mathscr A}_{\omega^n}(G)\to A$. Очевидно, тогда его можно продолжить до непрерывного гомоморфизма из ${\mathscr A}_{\omega^\infty}(G)$. Поскольку ${\mathscr A}_{\omega^\infty}(G)$ является пополнением ${\mathscr A}(G)$, такое продолжение единственно. Итак, универсальное свойство имеет место, и это доказывает теорему.
\end{proof}

Заметим, что непосредственно из определений следует, что существует непрерывный гомоморфизм $\chi\!:\widehat{\mathscr A}(G)\to {\widehat{{\mathscr A}}}(G)^{\mathrm{PI}}$ алгебр Аренса-Майкла (см. также теорему~\ref{redtog}).

\begin{pr}\label{PIH}
${\widehat{{\mathscr A}}}(G)$ и ${\widehat{{\mathscr A}}}(G)^{\mathrm{PI}}$  являются $\mathbin{\widehat{\otimes}}$-алгебрами Хопфа относительно операций, непрерывно продолженных с ${\mathscr A}(G)$, а отображение $\chi\!:\widehat{\mathscr A}(G)\to {\widehat{{\mathscr A}}}(G)^{\mathrm{PI}}$ является гомоморфизмом $\mathbin{\widehat{\otimes}}$-алгебр Хопфа.
\end{pr}
\begin{proof}
Согласно теореме~\ref{redtog} алгебра ${\mathscr A}_{\omega^\infty}(G)$  изоморфна ${\widehat{{\mathscr A}}}(G)$, в случае, если $\omega$ --- максимальный вес, и ${\widehat{{\mathscr A}}}(G)^{\mathrm{PI}}$  в случае, если $\omega$ --- максимальный вес с экспоненциальным искривлением на $N$. Из \cite[Lemma 1.9]{ArUAPI} следует, что в обоих случаях можно предполагать $\omega$ симметричным. Таким образом, выполнены условия предложения~\ref{ominftyAM2}, из которого следует, что ${\widehat{{\mathscr A}}}(G)$ и ${\widehat{{\mathscr A}}}(G)^{\mathrm{PI}}$ --- $\mathbin{\widehat{\otimes}}$-алгебры Хопфа, а отображения пополнения $\iota\!:{\mathscr A}(G)\to {\widehat{{\mathscr A}}}(G)$ и $\iota'\!:{\mathscr A}(G)\to {\widehat{{\mathscr A}}}(G)^{\mathrm{PI}}$ являются гомоморфизмами $\mathbin{\widehat{\otimes}}$-алгебр Хопфа (утверждения для ${\widehat{{\mathscr A}}}(G)$ также следуют из \cite[Proposition 6.7]{Pir_stbflat}).
Так как $\iota'=\chi\iota$, а образ $\iota$ плотен, то $\chi$ --- также гомоморфизм $\mathbin{\widehat{\otimes}}$-алгебр Хопфа.
\end{proof}

\section{Разложение в итерированные смэш-произведения алгебр Аренса-Майкла}
\label{s:DISPAM}

В этом разделе доказаны общие результаты о разложении в итерированное смэш-произведение и рассмотрены частные случаи --- $\widehat{{\mathscr A}}(G)$ и $\widehat{{\mathscr A}}(G)^{\mathrm{PI}}$.

\subsection*{Общий случай}

Начнём с общего утверждения. (Обозначения здесь и ниже --- те же, что в \S\,\ref{s:TDWED}.)

\begin{thm}\label{wtGDAF0}
Пусть $G$ --- связная линейная комплексная группа Ли, и $\omega_{max}$ --- субмультипликативный вес на~$G$, максимальный среди весов с экспоненциальным искривлением на некоторой нормальной интегральной подгруппе $N'$ такой, что $E\subset N'\subset N$. Тогда алгебра ${\mathscr A}_{\omega_{max}^\infty}(G)$ топологически изоморфна итерированному смэш-произведению:
$$
{\mathscr A}_{\omega_{max}^\infty}(G)\cong((\cdots ({\mathscr A}_{\omega_1^\infty}(F_1) \mathop{\widehat{\#}} {\mathscr A}_{\omega_2^\infty}(F_2))\mathop{\widehat{\#}}\cdots)
\mathop{\widehat{\#}} {\mathscr A}_{\omega_n^\infty}(F_n)),
$$
которое совпадает с итерированным смэш-произведением
$$
((\cdots ({\mathscr A}_{\omega_1^\infty}(F_1) \mathop{\#}\nolimits_{\mathsf{AM}} {\mathscr A}_{\omega_2^\infty}(F_2))\mathop{\#}\nolimits_{\mathsf{AM}}\cdots)
\mathop{\#}\nolimits_{\mathsf{AM}} {\mathscr A}_{\omega_n^\infty}(F_n))
$$
\emph{(}в обозначениях из теоремы~\ref{itsemdd0}\emph{)}.
При этом гомоморфизм ${{\mathscr A}}(G)\to {\mathscr A}_{\omega_{max}^\infty}(G)$ имеет вид, согласованный с этим разложением и разложением в следствии~\ref{GDAF}.
\end{thm}
\begin{proof}
Согласно теореме~\ref{itsemdd0} вес $\omega_{max}$ имеет итерированное разложение с точностью до эквивалентности. Более того, в силу \cite[Lemma 1.9]{ArUAPI}  можно предполагать, что $\omega_{max}$ симметричен. Тем самым его ограничение $\widetilde\omega_i$ на каждую из подгрупп~$G_i$, составляющих композиционный ряд, так же симметрично.  Веса $\omega_i$, участвующие в разложении, симметричнвы в силу их определения в~\eqref{omiz}. Итак, мы можем применить
теорему~\ref{smprde}, из которой следует существование изоморфизма ${\mathscr A}_{\omega_{max}^\infty}(G)$ и итерированного смэш-произведения. Соласованность с разложением в следствии~\ref{GDAF} следует из предложения~\ref{comdiaanco}.
\end{proof}

Теперь мы выпишем явный вид множителей  в теореме~\ref{wtGDAF0}.
Для этого понадобится следующая конструкция.

Для $s\ge 0$ положим
\begin{equation}
 \label{faAsdef}
\mathfrak{A}_s\!:=\Bigl\{a=\sum_{n=0}^\infty  a_n x^n\! :
\|a\|_{r,s}\!:=\sum_{n=0}^\infty |a_n|\frac{r^n}{n!^s}<\infty
\;\forall r>0\Bigr\},
\end{equation}
где $x$ --- формальная переменная.
Для нас важно, что $\mathfrak{A}_s$ является алгеброй Фреше-Аренса-Майкла  \cite[предложение~4]{ArRC} и, более того, $\mathbin{\widehat{\otimes}}$-алгеброй Хопфа \cite[Example 2.4]{AHHFG} (операции продолжаются по непрерывности с алгебры многочленов $\mathbb{C}[x]$).

\begin{lm}\label{1dimdsp}
\emph{(A)} Если $\omega(z)=1+|z|$, то ${\mathscr A}_{\omega^\infty}(\mathbb{C})$ топологически изоморфна $\mathbin{\widehat{\otimes}}$-алгебре Хопфа $\mathbb{C}[[x]]$ всех формальных степенных рядов от~$x$.

\emph{(B)} Если $\omega(z)=\exp(|z|^{1/s})$ для $s\in [1,\infty)$, то
${\mathscr A}_{\omega^\infty}(\mathbb{C})$ топологически изоморфна $\mathbin{\widehat{\otimes}}$-алгебре Хопфа $\mathfrak{A}_{s-1}$.

\emph{(C)} $\mathfrak{A}_0\cong \mathcal{O}(\mathbb{C})$.
\end{lm}
\begin{proof}
(A)~В силу \cite[Lemma 2.10]{ArAMN} для каждого $n\in\mathbb{N}$ дуальное пространство к банаховой алгебре ${{\mathscr A}}_{\omega^n}(\mathbb{C})$ изоморфно  банаховому пространству
\begin{equation}\label{cOom}
\mathcal{O}_{\omega^n}(\mathbb{C})\!:=\Bigl\{ f\in\mathcal{O}(\mathbb{C}) \!: |f|_{\omega^n}\!:=\sup_{z\in
\mathbb{C}}{\omega(z)}^{-n}{|f(z)|}<\infty\Bigr\}\quad(n\in\mathbb{N}).
\end{equation}

Так как $\omega(z)=1+|z|$, то, как легко видеть, $\mathcal{O}_{\omega^n}(\mathbb{C})$  является пространством многочленов степени не больше~$n$. Отсюда легко получить, что
${\mathscr A}_{\omega^n}(\mathbb{C})$ и $\mathbb{C}[x]/(x^{n+1})$ изоморфны как ассоциативные и, более того, как банаховы алгебы. Так как ${\mathscr A}_{\omega^\infty}(\mathbb{C})$ есть проективный предел последовательности ${\mathscr A}_{\omega^n}(\mathbb{C})$, то ${\mathscr A}_{\omega^\infty}(\mathbb{C})\cong \mathbb{C}[[x]]$ как локально выпуклое пространство. Непосредственная проверка показывает, что мы получили  изоморфизм $\mathbin{\widehat{\otimes}}$-алгебр Хопфа.

Доказательство части~(B) см. в~\cite[Lemma 4.5]{AHHFG}.
Часть~(C) следует непосредственно из~\eqref{faAsdef} (достаточно положить $s=0$).
\end{proof}

Теперь мы можем выписать разложение из теоремы~\ref{wtGDAF0} в явном виде.

\begin{thm}\label{exwtGDAF0}
Разложение из теоремы~\ref{wtGDAF0} может быть записано в следующем виде:
\begin{multline}\label{expfsmp}
 {\mathscr A}_{\omega_{max}^\infty}(G)\cong(\cdots (\mathbb{C}[[x_1]] \mathop{\widehat{\#}} \cdots \mathbb{C}[[x_p]])\mathop{\widehat{\#}}\mathfrak{A}_{m-1})\mathop{\widehat{\#}}\cdots\mathfrak{A}_{m-1})\mathop{\widehat{\#}}\cdots \\
 \cdots \mathop{\widehat{\#}}\mathfrak{A}_{1})\mathop{\widehat{\#}}\cdots\mathfrak{A}_{1})\mathop{\widehat{\#}}\mathcal{O}(\mathbb{C}))\mathop{\widehat{\#}}\cdots\mathcal{O}(\mathbb{C}))\mathop{\widehat{\#}} \widehat{{\mathscr A}}(L),
\end{multline}
где $p$ --- размерность $N'$, $m$ --- класс нильпотентности $B/N'$, а количество множителей вида $\mathfrak{A}_{j}$ для каждого $i$ равно размерности $j$-той факторалгебры нижнего центрального ряда~$B/N'$.
\end{thm}
\begin{proof}
Первые $p$ множителей в разложении изоморфны алгебре всех формальных степенных рядов ввиду того, что соответствующие веса имеют вид $z\mapsto 1+|z|$ (см.~\eqref{omiz} из теоремы~\ref{itsemdd0}) и части~(A) леммы~\ref{1dimdsp}.

Множители с номерами от $p+1$ до $n-1$ соответствуют весам вида  $z\mapsto \exp(|z|^{1/(j+1)})$, где $j$ пробегает значения от $m-1$ до $0$ в невозрастающем порядке с возможными повторениями  (см. также~\eqref{omiz}).
В силу части~(B) леммы~\ref{1dimdsp} получаем множители вида $\mathfrak{A}_{j}$, причём
количество повторений индекса~$j$ равно размерности $\mathfrak{q}_j/\mathfrak{q}_{j+1}$, где $\mathfrak{q}_1,\ldots,\mathfrak{q}_m$ --- нижний центральный ряд алгебры Ли группы $B/N'$ (см. замечание~\ref{seqwi}), а значит и размерности соотвествующей факторалгебры  нижнего центрального ряда самой~$B/N'$.

В заключение заметим, что $F_n=L$, где $L$ комплексно линейно редуктивна  и ${\mathscr A}_{\omega_n^\infty}(F_n)\cong\widehat{{\mathscr A}}(L)$, поскольку $\omega_n$ --- максимальный вес на~$L$.
\end{proof}

Таким образом, мы достигли поставленной цели --- представить ${\mathscr A}_{\omega_{max}^\infty}(G)$ в виде итерированного аналитического смэш-произведения с достаточно просто устроенными факторами --- алгебрами степенных рядов от одной переменной и  $\widehat{{\mathscr A}}(L)$, структура которой (также как и для  ${\mathscr A}(L)$) зависит только от теории представлений $L$, см. \cite[формула~(4.6)]{AHHFG}.

\subsection*{Два крайних случая}

Из теорем~\ref{wtGDAF0} и~\ref{redtog} получаем следующее утверждение.

\begin{thm}\label{wtGDAF}
Пусть $G$ --- связная линейная комплексная группа Ли. Тогда
алгебры $\widehat{{\mathscr A}}(G)$ и $\widehat{{\mathscr A}}(G)^{\mathrm{PI}}$ топологически изоморфны итерированному смэш-произведению вида
$$
((\cdots ({\mathscr A}_{\omega_1^\infty}(F_1) \mathop{\widehat{\#}} {\mathscr A}_{\omega_2^\infty}(F_2))\mathop{\widehat{\#}}\cdots)
\mathop{\widehat{\#}} {\mathscr A}_{\omega_n^\infty}(F_n)),
$$
которое совпадает с итерированным смэш-произведением
$$
((\cdots ({\mathscr A}_{\omega_1^\infty}(F_1) \mathop{\#}\nolimits_{\mathsf{AM}} {\mathscr A}_{\omega_2^\infty}(F_2))\mathop{\#}\nolimits_{\mathsf{AM}}\cdots)
\mathop{\#}\nolimits_{\mathsf{AM}} {\mathscr A}_{\omega_n^\infty}(F_n)),
$$
где $\omega_i$ определены в~\eqref{omiz}. В случае $\widehat{{\mathscr A}}(G)$ число $p$ равно размерности~$E$, а в случае  $\widehat{{\mathscr A}}(G)^{\mathrm{PI}}$ --- размерности~$N$.
Гомоморфизм ${{\mathscr A}}(G)\to \widehat{{\mathscr A}}(G)$ \emph{(}соответственно ${{\mathscr A}}(G)\to \widehat{{\mathscr A}}(G)^{\mathrm{PI}}$\emph{)} имеет вид, согласованный с этим разложением и разложением в следствии~\ref{GDAF}.
\end{thm}

Аналогично и теорема~\ref{exwtGDAF0} может быть записана в этих двух крайних случаях.

\begin{thm}\label{wtGDAFx}
Пусть $G$ --- связная линейная комплексная группа Ли.
Тогда $\widehat{{\mathscr A}}(G)$ топологически изоморфна итерированному смэш-произведению
\eqref{expfsmp}, где $p$ --- размерность~$E$ \emph{(}остальные обозначения те же, что и в теореме~\ref{exwtGDAF0}\emph{)}.
\end{thm}
\begin{proof}
В силу теоремы~\ref{redtog} случай, когда $N'=E$, соответствует оболочке Аренса-Майкла.   Тем самым $p$ равно размерности~$E$. Остальное ясно.
\end{proof}

\begin{thm}\label{wtPIGDAF}
Пусть $G$ --- связная линейная комплексная группа Ли.
Тогда ${\widehat{{\mathscr A}}}(G)^{\mathrm{PI}}$ топологически изоморфна итерированному смэш-произведению вида~\eqref{expfsmp},
где $p$ --- размерность~$N$, а $m=1$. А именно,
$$
{\widehat{{\mathscr A}}}(G)^{\mathrm{PI}}\cong(\cdots (\mathbb{C}[[x_1]] \mathop{\widehat{\#}} \cdots \mathbb{C}[[x_p]])\mathop{\widehat{\#}} \mathcal{O}(\mathbb{C}))\mathop{\widehat{\#}}\cdots\mathcal{O}(\mathbb{C}))\mathop{\widehat{\#}} \widehat{{\mathscr A}}(L).
$$
\end{thm}
\begin{proof}
В этом случае $N'=N$, см. теорему~\ref{redtog}. Заметим, что $B/N$ является абелевой, будучи нильпотентной и одновременно подгруппой комплексно линейно редуктивной группы Ли. Это означает, что её класс нильпотентности равен~$1$, а значит разложение из теоремы~\ref{exwtGDAF0} имеет требуемый вид.
\end{proof}

\begin{rem}
Утверждение теоремы~\ref{itsemdd0} в случае, когда $N'=N$,  также принимает более простой вид:  в части~(A) вместо последовательности разных весов $z\mapsto\exp(|z|^{1/w_{n-i}})$ получается последовательность одинаковых весов $z\mapsto\exp|z|$.  Это следует из замечания~\ref{seqwi}.
\end{rem}

\begin{rem}
В \cite{AHHFG} рассмотрен класс $\mathbin{\widehat{\otimes}}$-алгебр Хопфа, которые в то же время являются голоморфно конечно порождёнными алгебрами в смысле \cite{Pi14,Pi3}. В частности, показано, что $\widehat{\mathscr A}(G)$ является таковой. Хотя в полученном в теореме~\ref{wtGDAF} разложении в смэш-произведение на каждом шаге получается алгебра (Хопфа) Аренса-Майкла, тем не менее она не обязана быть голоморфно конечно порождённой.  Действительно, согласно \cite[предложение~13]{ArRC} алгебра $\mathfrak{A}_{t}$ не является таковой,  если $t\in(0,+\infty]$. Тем самым уже первый множитель в разложении может не принадлежать этому классу.

Вопрос от том, является ли ${\widehat{{\mathscr A}}}(G)^{\mathrm{PI}}$ голоморфно конечно порождённой, является более тонким. Можно предположить, что в общем случае ответ отрицательный.
\end{rem}

Сравнивая разложения в теоремах~\ref{wtGDAFx} и~\ref{wtPIGDAF}, получаем следующее утверждение.
\begin{co}
Пусть $G$ --- связная линейная комплексная группа Ли. Гомоморфизм  $\mathbin{\widehat{\otimes}}$-алгебр Хопфа
$\chi\!:{\widehat{{\mathscr A}}}(G)\to {\widehat{{\mathscr A}}}(G)^{\mathrm{PI}}$ является изоморфизмом тогда и только тогда, когда $E=N$.
\end{co}

\begin{rem}
Утверждение следствия можно переформулировать иначе: \emph{все пополнения ${\mathscr A}(G)$ относительно всевозможных непрерывных субмультипликативных преднорм являются \emph{PI}-алгебрами тогда и только тогда, когда $E=N$.} Формулируя этот критерий, автор исполняет обещание, данное в конце~\cite{ArPiLie}.
\end{rem}

В заключение рассмотрим два простых примера, демонстрирующие случаи, когда $\chi\!:{\widehat{{\mathscr A}}}(G)\to {\widehat{{\mathscr A}}}(G)^{\mathrm{PI}}$ является топологическим изоморфизмом и когда нет.

\begin{exm}
Пусть $\mathfrak{g}$ --- двумерная комплексная алгебра Ли с базисом $\{e_1, e_2\}$ и соотношением $[e_1, e_2] = e_2$.
Так как она разрешима, то $\mathfrak{n}=[\mathfrak{g}, \mathfrak{g}] = \mathbb{C} e_2$. С другой стороны, нетрудно проверить, что $\mathfrak{e}$ (алгебра Ли
группы~$E$) совпадает с $[\mathfrak{g}, \mathfrak{g}]$, см. \cite[Example 5.13]{ArAnF}. Пусть $G$ --- однозвязная группа Ли, такая что $\mathfrak{g}$ --- её алгебра Ли.
Так как $\mathfrak{e}=\mathfrak{n}$, то ${\widehat{{\mathscr A}}}(G)\to {\widehat{{\mathscr A}}}(G)^{\mathrm{PI}}$ --- топологический изоморфизм. Обе алгебры изоморфны
$\mathbb{C}[[e_2]]\mathop{\widehat{\#}}\mathcal{O}(\mathbb{C})$.
\end{exm}

\begin{exm}
Пусть $\mathfrak{g}$ --- трёхмерная алгебра Гейзенберга, т.е. порождена базисом $\{e_1, e_2, e_3\}$ и соотношениями $[e_1, e_2] = e_3$, $[e_1, e_2] = 0$, $[e_1, e_3] = 0$.  Пусть $G$ --- трёхмерная группа Гейзенберга, т.е. однозвязная группа Ли, такая что $\mathfrak{g}$ --- её алгебра Ли. Заметим, что из нильпотентности $G$ следует, что её экспоненциальный радикал тривиален, так же как и его алгебра Ли $\mathfrak{e}$. С другой стороны, $\mathfrak{n}=[\mathfrak{g}, \mathfrak{g}]=\mathbb{C} e_3$. Так как $\mathfrak{e}\ne\mathfrak{n}$, то ${\widehat{{\mathscr A}}}(G)\to {\widehat{{\mathscr A}}}(G)^{\mathrm{PI}}$ не является топологическим изоморфизмом. Разложения имеют следующий вид:
$$
{\widehat{{\mathscr A}}}(G)\cong (\mathfrak{A}_1\mathop{\widehat{\#}}\mathcal{O}(\mathbb{C}))\mathop{\widehat{\#}} \mathcal{O}(\mathbb{C})\quad\text{и}\quad{\widehat{{\mathscr A}}}(G)^{\mathrm{PI}}\cong(\mathbb{C}[[e_3]]\mathop{\widehat{\#}}\mathcal{O}(\mathbb{C}))\mathop{\widehat{\#}} \mathcal{O}(\mathbb{C}).
$$
Заметим, что вторая алгебра совпадает с алгеброй ``формально-радикальных функций'', рассмотренных Доси в  \cite{Do09C,DoSb,Do10A,Do10B}. Такое совпадение не случайно --- можно показать, что  ${\widehat{{\mathscr A}}}(G)^{\mathrm{PI}}$ изоморфна алгебре формально-радикальных функций для любой нильпотентной $\mathfrak{g}$. Это требует дополнительных рассуждений и будет опубликовано отдельно.
\end{exm}

\end{document}